  \documentclass[12pt]{amsart}
  \usepackage{latexsym}
  \usepackage[all]{xy}
  \usepackage{amsfonts}
  \usepackage{amsthm}
  \usepackage{amsmath}
  \usepackage{amssymb}
  \usepackage[mathscr]{eucal}
  \numberwithin{equation}{section}

  \def\sw#1{{\sb{(#1)}}}
  
  \def\su#1{{\sp{[#1]}}} %

  \def\<{{\langle}}
  \def\>{{\rangle}}
  
  \def\eps{\varepsilon}

  \def\note#1{{}}
  \def\can{{\rm can}}
  
  \def\note#1{}

  \def\cM{{\mathcal M}}

  \def\cC{{\mathfrak C}}
  
  \def\cD{{\mathfrak D}}
   \def\cE{{\mathfrak E}}

  \def\lhom#1#2#3{{{\rm Hom}\sb{#1-}(#2,#3)}}
  \def\rhom#1#2#3{{{\rm Hom}\sb{-#1}(#2,#3)}}
\def\hom#1#2#3{{{\rm Hom}\sb{#1}(#2,#3)}}
  
  \def\rend#1#2{{{\rm End}\sb{-#1}(#2)}}

  \def\Rend#1#2{{{\rm End}\sp{-#1}(#2)}}
  
  \def\Rhom#1#2#3{{{\rm Hom}\sp{-#1}(#2,#3)}}

  \def\LRhom#1#2#3#4{{{\rm Hom}\sp{#1,#2}(#3,#4)}}

  \def\can{{\rm can}}

  \def\cor{{\bf Crg}}

  \def\beq{\begin{equation}}
  \def\eeq{\end{equation}}
  \def\DC{{\Delta_\cC}}
  \def \eC{{\eps_\cC}}
  \def\DD{{\Delta_\cD}}
  \def \eD{{\eps_\cD}}

  \def\DSD{{\Delta_{\Sigma[\cD]}}}
  \def \eSD{{\eps_{\Sigma[\cD]}}}

  \def\Desc{{\bf Desc}}
  \def\Rep{{\bf Rep}}
  \def\MRep{{\bf Rep^M}}
  \def\id{{I}}

  \def\ut{{\otimes}}
  \def\ot{{\otimes}}
  
  \def\Hom{\mbox{\rm Hom}\,}
  \def\roM{\varrho^{M}}
\def\Nro{{}^{N}\!\varrho}
  
  \newcommand{\Ra}{\Rightarrow}
  \def\bSi{\mathbf{\Sigma}}
    \def\tbSi{\mathbf{\widetilde{\Sigma}}}
    \def\tSi{\widetilde{\Sigma}}
    \def\tfks{\widetilde{\mathfrak{s}}}
    \def\tfkw{\widetilde{\mathfrak{w}}}
\def\tsi{\tilde{\sigma}}
    
    \def\hSi{\widehat{\Sigma}}

\def\te{\tilde{e}}

\def\ts{\tilde{s}}

 \def\bXi{\mathbf{\Xi}}
    \def\tbXi{\mathbf{\widetilde{\Xi}}}
    \def\tXi{\widetilde{\Xi}}

 \newcommand{\rmod}[1]{\mathcal{M}_{#1}}

\newcommand{\rcomod}[1]{\mathcal{M}^{#1}}
\newcommand{\lcomod}[1]{{}^{#1}\mathcal{M}}
\newcommand{\bimod}[2]{{}_{#1}\mathcal{M}_{#2}}

\newcommand{\coring}[1]{\mathfrak{#1}}
\newcommand{\tensor}[1]{\otimes_{#1}}

\newcommand{\rcomatrix}[2]{#2^* \tensor{#1} #2}
\newcommand{\bara}[1]{\overline{#1}}
\newcommand{\cotensor}[1]{\square_{#1}}
\newcommand{\REM}{\mathsf{REM}(\mathsf{Bim})}
\newcommand{\LEM}{\mathsf{LEM}(\mathsf{Bim})}
\newcommand{\fREM}{\mathsf{fREM}(\mathsf{Bim})}
\newcommand{\fLEM}{\mathsf{fLEM}(\mathsf{Bim})}
\newcommand{\fk}[1]{\mathfrak{#1}}
\newcommand{\Bim}{\mathsf{Bim}}

  \newcounter{zlist}
  \newenvironment{zlist}{\begin{list}{(\arabic{zlist})}{
  \usecounter{zlist}\leftmargin2.5em\labelwidth2em\labelsep0.5em
  \topsep0.6ex
  \parsep0.3ex plus0.2ex minus0.1ex}}{\end{list}}

  \newcounter{blist}
  \newenvironment{blist}{\begin{list}{(\alph{blist})}{
  \usecounter{blist}\leftmargin2.5em\labelwidth2em\labelsep0.5em
  \topsep0.6ex 
  \parsep0.3ex plus0.2ex minus0.1ex}}{\end{list}}

  \newcounter{rlist}

  \def\Label#1{\label{#1}\ifmmode\llap{[#1] }\else
  \marginpar{\smash{\hbox{\tiny [#1]}}}\fi}
  \def\Label{\label}

  \newtheorem{proposition}{Proposition}[section]
  \newtheorem{lemma}[proposition]{Lemma}
  
  \newtheorem{corollary}[proposition]{Corollary}
  \newtheorem{theorem}[proposition]{Theorem}

  \theoremstyle{definition}
  \newtheorem{definition}[proposition]{Definition}
  \newtheorem{example}[proposition]{Example}

  \theoremstyle{remark}
  \newtheorem{remark}[proposition]{Remark}

\newcounter{c}
\renewcommand{\[}{\setcounter{c}{1}$$}
 \newcommand{\etyk}[1]{\vspace{-7.4mm}$$\begin{equation}\Label{#1}
\addtocounter{c}{1}}
\renewcommand{\]}{\ifnum \value{c}=1 $$\else \end{equation}\fi}
\begin{document}
\baselineskip 19pt

  \title{The bicategories of corings}
  \author{Tomasz Brzezi\'nski}
  \address{ Department of Mathematics, University of Wales Swansea,
  Singleton Park, \newline\indent  Swansea SA2 8PP, U.K.}
  \email{T.Brzezinski@swansea.ac.uk}
  \urladdr{http//www-maths.swan.ac.uk/staff/tb}
 \author{ L.\ El Kaoutit}
 \address{Departamento de {\'A}lgebra,  Universidad de
Granada, E18071 Granada, Spain}
\email{kaoutit@ugr.es}
 \author{J.\ G\'omez-Torrecillas}
 \address{Departamento de {\'A}lgebra,  Universidad de
Granada, E18071 Granada, Spain} \email{gomezj@ugr.es}
\urladdr{http://www.ugr.es/local/gomezj}

\begin{abstract}
To a $B$-coring and a $(B,A)$-bimodule that is finitely generated
and projective as a right $A$-module an $A$-coring is associated.
This new coring is termed a {\em base ring extension of a coring
by a module}. We study how the properties of a bimodule such as
separability and the Frobenius properties are reflected in the
induced base ring extension coring. Any bimodule that is finitely
generated and projective on one side, together with a map of
corings over the same base ring, lead to the notion of a {\em
module-morphism}, which extends the notion of a morphism of
corings (over
 different base rings). A module-morphism of corings induces functors between the categories of comodules.
 These functors are termed {\em pull-back} and {\em push-out} functors respectively
 and thus relate categories of comodules of different corings.
We study when the pull-back
functor is fully faithful and when it is an equivalence. A
 generalised descent associated to a
 morphism of corings is introduced.
 We define a {\em category of
 module-morphisms}, and show that push-out functors are naturally isomorphic
 to each other if and only if the corresponding module-morphisms are mutually
isomorphic. All these topics are studied within a unifying language of
{\em bicategories} and the extensive use is made of interpretation of corings
as comonads in the bicategory $\Bim$ of bimodules and
module-morphisms as 1-cells in the associated bicategories of comonads
in $\Bim$.
\end{abstract}
  \maketitle
{\small \tableofcontents }
  \section{Introduction}
  \subsection{Motivation and overview}
The aim of this paper is to study properties of corings and
functors between comodule categories from a {\em bicategorical}
point of view, and thus argue that bicategories provide a natural
and unifying point of view on corings.

In context of corings, bicategories arise in a very natural way.
The categorical information about rings is contained in a
bicategory of bimodules $\Bim$ in which objects (0-cells) are
rings, 1-cells are bimodules (with composition given by the tensor
product) and 2-cells are bilinear maps. This is the most
fundamental example of a bicategory, and provides the ideal set-up
for studying problems such as Morita theory. Building upon a
pioneering work of Street \cite{Str:for}, Lack and Street
\cite{Lack/Street:2002} have considered the bicategory
$\mathsf{EM}(\mathsf{B})$ obtained as
the
free completion under Eilenberg-Moore objects of
a given  bicategory $\mathsf{B}$.
When $\mathsf{B}$ is taken to be $\Bim$, the suitable dual
of $\mathsf{EM}(\mathsf{B})$ is  a bicategory
in which objects are corings. The resulting bicategory (and its
sub-bicategories) is the main object of studies in the present
paper.

Our motivation for studying the bicategory of corings is twofold,
deeply rooted in  non-commutative geometry. First, it appears that
there is a growing appreciation for the language of bicategories
in non-commutative geometry. For example, in a recent paper
\cite{Man:rea} Manin argues that the classification of vector
bundles over the non-commutative torus (or the K-theory of this
torus) is best explained in terms of a Morita bicategory
associated with this torus. This particular example of the role
that bicategories play in non-commutative geometry can be seen in
a much wider context, as bicategories appear very naturally in
quantisation of Poisson manifolds in terms
 of $C^*$-algebras or in the theory of von Neumann algebras \cite{Lan:bic}.

  Second motivation originates from the appearance of corings in non-commutative
  algebraic geometry. Here corings feature in two different ways. On one hand, if an $A$-coring
  $\cC$ is flat as a left $A$-module, then the category of its right comodules
  is a Grothendieck category, hence a non-commutative space or a
  non-commutative quasi-scheme in the
  sense of Van den Bergh \cite{Ber:blo} and Rosenberg \cite{Ros:alg} \cite{Ros:sch}.
  Natural isomorphism classes of functors between Grothendieck categories play
  the role of maps between non-commutative spaces (cf.\ \cite{Smi:sub}). In relation
  to corings, one needs to consider bimodules between corings that induce functors
  between corresponding comodule categories. Natural maps between these functors
  arise from morphisms between corresponding bimodules. It turns out that to study all
  these structures in a unified way one is led to considering a suitable bicategory. In another
  approach to non-commutative algebraic geometry, certain classes of corings
  appear as {\em covers} of non-commutative spaces \cite{KonRos:smo}. Bimodules
  between corings can then be understood as a change of cover of the underlying
  space. The change of cover affects the corresponding quasi-scheme, i.e.\ the
  category of comodules. Again to study these effects in a uniform way, one
  should study a bicategory in which corings are 0-cells.

Our presentation and the choice of topics for the present paper
are motivated by the above geometric interpretation. We begin by
applying (the comonadic version of) the Lack and Street
construction to the bicategory of bimodules and describe
explicitly the resulting bicategory of corings $\REM$. With an eye
on the interpretation of corings as covers of non-commutative
spaces we introduce the bicategory $\fREM$, by restricting 1-cells
in $\REM$ to those that arise  from adjoint pairs in $\Bim$
(finitely generated and projective modules). The Lack and Street
construction has an obvious `left-sided' version resulting in
bicategories $\LEM$ and $\fLEM$.
We show that there is a duality between the
hom-categories of $\fLEM$ and $\fREM$.

Next we introduce the notion of a base extension of a coring by a
bimodule (a `change of cover' of a non-commutative space). We show
that any such base extension gives rise to 1-cells in $\fLEM$ and
$\fREM$. We study in what way  properties of a bimodule such as
separabality and the Frobenius property are reflected by the
resulting base extension coring.

The appearance of  1-cells in $\fLEM$ and $\fREM$ associated to a
base extension of a coring by a bimodule leads to the notion of  a
{\em module-morphism} of corings as a pair consisting of  a
bimodule that is finitely generated and projective as a right
module, and of a coring map.  We introduce two functors between
categories of comodules induced by a  module-morphism. Following
the geometric interpretation these are called  a {\em push-out}
and a {\em pull-back} functors. In Section~\ref{sec.equiv} we
determine when these  functors are inverse equivalences, and also
we define and give basic properties of a generalised descent
associated to a morphism of corings.

\subsection{Notation and preliminaries}
We work over a commutative ring $k$ with a unit. All algebras are
over $k$, associative and with a unit. The identity morphism for
an object $V$ is also denoted by $V$. For a ring (algebra) $R$,
the category of right $R$-modules and right $R$-linear maps is
denoted by $\cM_R$. Symmetric notation is used for left modules.
As is customary, we often write $M_R$ to indicate that $M$ is a
right $R$-module, etc.
The dual module of $M_R$, consisting of all
$R$--linear maps from $M_R$ to $R_R$, is denoted by $M^*$, while
the dual of ${}_RN$ is denoted by ${}^*N$. The multiplication in
the endomorphism ring of a right module (comodule) is given by
composition of maps, while the multiplication in the endomorphism
ring of a left module (comodule) is given by opposite composition
(we always write argument to the right of a function). The symbol
$-\ot_R-$ between maps and modules denotes the tensor product
bifunctor over the algebra $R$.

Let $A$ be an algebra. The comultiplication of an $A$-coring $\cC$
is denoted by  $\DC :\cC\to \cC\ot_A\cC$, and its counit by
$\eC:\cC\to A$. To indicate the action of $\DC$ on elements we use
the Sweedler sigma notation, i.e.\ for all $c\in \cC$,
$$
\DC(c) = \sum c\sw1\ot_A c\sw 2, \quad (\DC\ot_A\cC)\circ\DC(c) (\cC\ot_A\DC)\circ\DC(c) = \sum c\sw1\ot_A c\sw 2\ot_A c\sw 3,
$$
etc. Gothic capital letters always denote corings. To indicate
that $\cC$ is an $A$-coring we often write $(\cC\! : \!A)$. The
category of right $\cC$-comodules and right $\cC$-colinear maps is
denoted by $\cM^\cC$. Recall that $\cM^\cC$ is built upon the
category of right $A$-modules, in the sense that there is a
forgetful functor $\cM^\cC\to \cM_A$. In particular, any right
$\cC$-comodule is also a right $A$-module, and any right
$\cC$-comodule map is right $A$-linear. For a right $\cC$-comodule
$M$, $\varrho^M:M\to M\ot_A\cC$ denotes the coaction, and $\Rhom
\cC MN$ is the $k$-module of $\cC$-colinear maps $M\to N$. On
elements, the action of $\varrho^M$ is expressed by Sweedler's
sigma-notation $\varrho^M(m) = \sum m\sw 0\ot m\sw 1$. Symmetric
notation is used for left $\cC$-comodules. In particular, the
coaction of a left $\cC$-comodule $N$ is denoted by ${}^N\!\varrho$,
and, on elements, by ${}^N\!\varrho(n) = \sum n\sw{-1}\ot n\sw 0\in
\cC\ot_A N$.

For any $A$-coring $\cC$, the dual module $\cC^*= \rhom A\cC A$ is
a $k$-algebra with the product $f*g(c) = \sum f(g(c\sw 1)c\sw 2)$
and unit $\eC$. This is known as a {\em right dual ring of $\cC$}.
Similarly, the dual module ${}^*\cC= \lhom A\cC A$ is a
$k$-algebra with the product $f*g(c) = \sum f(c\sw 1g(c\sw 2))$
and unit $\eC$. This is known as a {\em left dual ring of $\cC$}.
The $k$-linear map $\iota_A: A\to {}^*\cC$,
$a\mapsto[c\mapsto\eps_\cC(ca)]$ is an anti-algebra map.

Given $A$-corings $\cC$, $\cD$, a {\em morphism of $A$-corings}
$\cD\to \cC$ is an $A$-bimodule map $\gamma: \cD\to\cC$ such that
$\DC\circ\gamma = (\gamma\ot_A\gamma)\circ\DD$ and $\eC\circ\gamma = \eD$.
The category of $A$-corings is denoted by $A$-\cor

If $\cD$ is a $B$-coring and $\alpha:B\to A$ is an algebra map,
then one views $A$ as a $B$-bimodule via $\alpha$ and defines an
$A$-coring structure on the $A$-bimodule $A_\alpha[\cD] := A\ot_B \cD\ot_B A$ by
\begin{eqnarray*}
\Delta_{A_\alpha[\cD]}: A_\alpha[\cD]&\to& A_\alpha[\cD]\ot_AA_\alpha[\cD]
\simeq A\ot_B \cD\ot_B A\ot_B \cD\ot_B A, \\
 a\ot_B d\ot_B a'&\mapsto& \sum a\ot_B d\sw 1\ot_B 1\ot_B d\sw 2\ot_B a',
\end{eqnarray*}
and $\eps_{A_\alpha[\cD]}: a\ot d\ot a'\mapsto a\eD(d)a'$. $A_\alpha[\cD]$
is known as a {\em base ring extension of $\cD$}. The construction of a base
ring extension allows one to consider morphisms of corings over different rings.
Given corings $(\cC\! : \!A)$ and $(\cD\! : \!B)$ a {\em morphism of corings}
$(\cD\! : \!B)\to (\cC\! : \!A)$ is a pair $(\gamma,\alpha)$,
where $\alpha :B\to A$ is an algebra map and  $\gamma: \cD\to \cC$ is a
$B$-bimodule map such that the induced map $\tilde{\gamma}:A_\alpha[\cD] \to \cC$,
$a\ot_Bd\ot_Ba'\mapsto a\gamma(d)a'$ is a morphism of $A$-corings.

Recall that,
given a right $\cC$-comodule $M$ and a left $\cC$-comodule $N$ one defines a
{\em cotensor product} $M\Box_\cC N$ by the following exact  sequence of $k$-modules:
$$ \xymatrix{0\ar[r]& M\Box_\cC N\ar[r]& M\ot_A N\ar[r]^{\omega_{M,N}\quad}&
            M\ot_A  \cC\ot_A  N,}$$
where $\omega_{M,N} =\roM \ut_A N - M
\ut_A \Nro$, and  $\roM$ and $\Nro$ are coactions.
Suppose that $\cC$ is flat as a left $A$-module. A left
$\cC$-comodule $N$ is said to be {\em coflat} (resp.\ {\em
faithfully coflat}), if the cotensor functor $-\Box_\cC N:
\mathcal{M}^\cC\to \mathcal{M}_k$ preserves (resp.\ preserves and
reflects) exact sequences in $\mathcal{M}^\cC$.
A detailed account of the theory of corings and comodules can be
found in \cite{BrzWis:cor}.

Henceforth, $\{e_i,e_i^*\}$, with $e_i\in \Sigma$, $e^*_i\in
\Sigma^*$ always denotes a finite dual bases for a finitely
generated and projective module $\Sigma_A$.

\section{The bicategories of corings}\label{sec.bicat}

For general definitions of bicategories and their morphisms we
refer the reader to the fundamental paper \cite{Benabou:1967}.
Following \cite[(2.5)]{Benabou:1967} and the conventions adopted
there, the bicategory of bimodules $\Bim$ is defined as follows.
Objects (i.e. $0$-cells) are algebras $A$, $B$,...,  $1$-cells
from $A$ to $B$ are objects of the category $\bimod{B}{A}$ of
$(B,A)$-bimodules, and $2$-cells are  bilinear maps. The
composition of $1$-cells is given by the tensor product of
bimodules, and the identity $1$-cell of $A$ is the regular
bimodule ${}_{A}A_A$. Given a bicategory $\mathsf{B}$, the
transpose bicategory of $\mathsf{B}$ is the bicategory
$\mathsf{B}^{op}$ obtained from $\mathsf{B}$ by reversing
$1$-cells, while the conjugate bicategory $\mathsf{B}^{co}$ is
obtained by reversing $2$-cells, and the 
bicategory
$\mathsf{B}^{coop}$ is obtained by reversing both (cf. \cite[\S
3]{Benabou:1967}).

\subsection{The right bicategory of corings}

Following \cite[p.\ 249]{Lack/Street:2002} (cf.\ \cite{Str:for}),
the  bicategory $\REM \,:=\,\mathsf{EM}\left(
\Bim^{coop}\right)^{coop}$  consists of
the following data:

$\bullet$ \emph{Objects}: Corings $(\cC\! :\! A)$ (i.e.\ $\cC$ is an
$A$-coring).

$\bullet$ $1$-\emph{cells}: A $1$-cell from $(\cC\! :\! A)$ to
$(\cD\! :\! B)$ is a pair $(\Sigma,\fk{s})$ consisting of a
$(B,A)$-bimodule $\Sigma$ and a $(B,A)$-bilinear map $\fk{s}:
\cD\tensor{B}\Sigma \rightarrow \Sigma\tensor{A}\cC$ rendering
commutative the following diagrams
$$\xymatrix@R=30pt@C=30pt{\cD\tensor{B}\Sigma \ar@{->}^-{\fk{s}}[r]
\ar@{->}|-{\varepsilon_{\cD}\tensor{B}\Sigma}[d] &
\Sigma\tensor{A}\cC\ar@{->}|-{\Sigma\tensor{A}\varepsilon_{\cC}}[d] \\
B\tensor{B}\Sigma \ar@{->}^-{\simeq}[r] & \Sigma\tensor{A}A,} \quad
\xymatrix@R=30pt@C=40pt{\cD\tensor{B}\Sigma \ar@{->}^-{\fk{s}}[r]
\ar@{->}|-{\Delta_{\cD}\tensor{B}\Sigma}[d] & \Sigma\tensor{A}\cC
\ar@{->}^-{\Sigma\tensor{A}\Delta_{\cC}}[r] &
\Sigma\tensor{A}\cC\tensor{A}\cC
\\ \cD\tensor{B}\cD\tensor{B} \Sigma \ar@{->}^-{\cD\tensor{B}\fk{s}}[r] &
\cD\tensor{B}\Sigma\tensor{A} \cC
\ar@{->}_-{\fk{s}\tensor{A}\cC}[ru] & }$$ The identity $1$-cell of
an object $(\cC\! :\! A)$ is given by $(A,\cC)$.

$\bullet$ $2$-\emph{cells}: Given $1$-cells $(\Sigma,\fk{s})$ and
$(\tSi,\tfks)$  from $(\cC\! :\! A)$ to $(\cD\! :\! B)$,
$2$-cells are defined as $(B,A)$-bilinear maps $\fk{a}:
\cD\tensor{B}\Sigma \to \tSi$  rendering commutative the
following diagram
$$\xymatrix@R=40pt@C=70pt{ \cD\tensor{B}\Sigma
\ar@{->}^-{\Delta_{\cD}\tensor{B}\Sigma}[r]
\ar@{->}|-{\Delta_{\cD}\tensor{B}\Sigma}[d] &
\cD\tensor{B}\cD\tensor{B}\Sigma
\ar@{->}^-{\cD\tensor{B}\fk{s}}[r] &
\cD\tensor{B}\Sigma\tensor{A}\cC \ar@{->}|-{\fk{a}\tensor{A}\cC}[d] \\
\cD\tensor{B}\cD\tensor{B}\Sigma
\ar@{->}^-{\cD\tensor{B}\fk{a}}[r] & \cD\tensor{B}\tSi
\ar@{->}^-{\tfks}[r] & \tSi\tensor{A}\cC. }$$

The category consisting of all $1$ and $2$-cells from $(\cC\! :\! A)$
to $(\cD\! :\! B)$ is denoted by
${}_{(\cD : B)}\mathcal{R}_{(\cC : A)}$. The composition of $2$-cells
is defined as follows: Let $(\Sigma,\fk{s})$, $(\tSi,\tfks)$ be
$1$-cells from $(\cC\! :\! A)$ to $(\cD\! :\! B)$, and $(W,\fk{w})$,
$(\widetilde{W},\tfkw)$ be $1$-cells from $(\cE\! :\! C)$ to $(\cC\! :\! A)$.
The composition of $1$-cells leads to the following $1$-cells
from $(\cE\! :\! C)$ to $(\cD\! :\! B)$:
$$\left(\,\Sigma\tensor{A}W, (\Sigma\tensor{A}\fk{w}) \circ
(\fk{s}\tensor{A}W)\,\right) \text{ and }
\left(\,\tSi\tensor{A}\widetilde{W},(\tSi\tensor{A}\tfkw) \circ
(\tfks\tensor{A}\widetilde{W})\,\right).
$$
If $\fk{a}:
\cD\tensor{B}\Sigma \to \tSi$ and $\fk{b}: \cC\tensor{A}W \to
\widetilde{W}$ are $2$-cells, then the horizontal composition $\fk{a}\otimes \fk{b}$ is
given by
$$\xymatrix@R=40pt@C=70pt{ \cD\tensor{B}\Sigma\tensor{A}W
\ar@{->}^-{\Delta_{\cD}\tensor{B}\Sigma\tensor{A}W}[r]
\ar@{-->}_{\fk{a}\tensor{}\fk{b}}[rrdd] &
\cD\tensor{B}\cD\tensor{B}\Sigma\tensor{A}W
\ar@{->}^-{\cD\tensor{B}\fk{a}\tensor{A}W}[r] & \cD\tensor{B}\tSi\tensor{A}W \ar@{->}|-{\tfks\tensor{A}W}[d] \\
& & \tSi\tensor{A}\cC \tensor{A} W \ar@{->}|-{\tSi\tensor{A}\fk{b}}[d] \\
& & \tSi\tensor{A}\widetilde{W}  }$$

Every $1$-cell in $\mathsf{REM (Bim)}$ defines a functor between
categories of right comodules. This statement is made explicit in the following

\begin{proposition}\label{push-Out}
Let $(\Sigma,\fk{s})$ be a $1$-cell from $(\cC\! :\! A)$ to $(\cD\! :\! B)$ in
the bicategory $\REM$. There is a functor $\bSi_{\circ}:
\rcomod{\cD} \longrightarrow \rcomod{\cC}$ sending $$\left(\,
(M,\varrho^M) \mapsto
(M\tensor{B}\Sigma,\varrho^{M\tensor{B}\Sigma}=(M\tensor{B}\fk{s})
\circ (\varrho^M\tensor{B}\Sigma) \,\right), \quad \left( f
\mapsto f\tensor{B}\Sigma \right).$$ In particular
$\cD\tensor{B}\Sigma$ admits a structure of a
$(\cD,\cC)$-bicomodule.
\end{proposition}
\begin{proof}
For any object $(\cC\! :\! A)$, the category
${}_{(k : k)}\mathcal{R}_{(\cC : A)}$ is isomorphic to the category $\rcomod{\cC}$ of
right $\cC$-comodules. The functor $\bSi_{\circ}$ is then identified with the horizontal composition functor
$-\otimes(\Sigma,\fk{s}): {}_{(k : k)}\mathcal{R}_{(\cD : B)}
\longrightarrow {}_{(k : k)}\mathcal{R}_{(\cC : A)}.$ View $\cD\tensor{B}\Sigma$ as a left $\cD$-comodule with  the coaction ${}^{\cD\ot_B\Sigma}\varrho=\Delta_{\cD}\ot_B\Sigma$, and compute
\begin{eqnarray*}
  (\cD\tensor{B}\varrho^{\cD\tensor{B}\Sigma}) \circ (\Delta_{\cD}\tensor{B}\Sigma)
  &=& (\cD\tensor{B}\cD\tensor{B}\fk{s}) \circ
  (\cD\tensor{B}\Delta_{\cD}\tensor{B}\Sigma) \circ (\Delta_{\cD}\tensor{B}\Sigma) \\
   &=& (\cD\tensor{B}\cD\tensor{B}\fk{s}) \circ
  (\Delta_{\cD}\tensor{B}\cD\tensor{B}\Sigma) \circ (\Delta_{\cD}\tensor{B}\Sigma) \\
   &=& (\Delta_{\cD}\tensor{B}\Sigma\tensor{A}\cC) \circ
   (\cD\tensor{B}\fk{s}) \circ (\Delta_{\cD}\tensor{B}\Sigma) \\
   &=& (\Delta_{\cD}\tensor{B}\Sigma\tensor{A}\cC) \circ
   \varrho^{\cD\tensor{B}\Sigma}.
\end{eqnarray*}
Hence ${}^{\cD\ot_B\Sigma}\varrho$ is right $\cC$-colinear, i.e.\
$\cD\tensor{B}\Sigma$ is a $(\cD,\cC)$-bicomodule, as stated.
\end{proof}

As observed in \cite[p.\ 249]{Lack/Street:2002}, 2-cells in a
bicategory of monads can be defined in a {\em reduced} or an {\em
unreduced} form. This bijective correspondence can be dualised to
bicategories of comonads, and, in the case of $\REM$, can be
interpreted in terms of bicomodules.

\begin{proposition}\label{red-unred}
Let $(\Sigma,\fk{s})$ and $(\tSi,\tfks)$ be two $1$-cells from
$(\cC\! :\! A)$ to $(\cD\! :\! B)$ in the bicategory $\REM$. There is a
bijection
$$
{}_{(\cD : B)}\mathcal{R}_{(\cC : A)}\left({(\Sigma,\fk{s})},{(\tSi,\tfks)}\right)
\,\, \simeq \,\,
\Hom^{\cD,\cC}(\cD\tensor{B}\Sigma,\cD\tensor{B}\tSi),
$$
 where $\cD\tensor{B}\Sigma$ and $\cD\tensor{B}\tSi$ are  $(\cD,\cC)$-bicomodules by Proposition~\ref{push-Out}. Explicitly,
$$
 \left( \xymatrix{  \fk{a} \ar@{|->}[r] &
(\cD\tensor{B}\fk{a}) \circ (\Delta_{\cD}\tensor{B}
\Sigma)}\right), \quad \left( \xymatrix{
(\varepsilon_{\cD}\tensor{B}\tSi) \circ f & \ar@{|->}[l] f }
\right).
$$
\end{proposition}
\begin{proof}
We only need to prove that the mutually inverse maps are well
defined. For a $2$-cell $\fk{a}: \cD\tensor{B}\Sigma \to
\tSi$,  the left
   $\cD$-colinearity of $((\cD\tensor{B}\fk{a}) \circ
   (\Delta_{\cD}\tensor{B}\Sigma))$ follows by the following simple calculation that uses the coassociativity of $\DD$:
\begin{eqnarray*}
   (\Delta_{\cD}\tensor{B}\tSi) \!\circ\! (\cD\tensor{B}\fk{a}) \!\circ\!(\Delta_{\cD}\tensor{B}\Sigma)\!\!\!
   &=&\!\!\!
   (\cD\tensor{B}\cD\tensor{B}\fk{a}) \!\circ\! (\Delta_{\cD}\tensor{B}\cD\tensor{B}\Sigma)
   \!\circ\! (\Delta_{\cD}\tensor{B}\Sigma)  \\
      &=&\!\!\! (\cD\tensor{B}( \,(\cD\tensor{B}\fk{a})  \!\circ\!(\Delta_{\cD}\tensor{B}\Sigma)\,))
      \!\circ\! (\Delta_{\cD}\tensor{B}\Sigma).
\end{eqnarray*}
Using the above calculation (to derive the second equality) and the fact that $\fk{a}$ is a 2-cell (to derive the third equality), one computes
\begin{eqnarray*}
  \varrho^{\cD\tensor{B}\tSi} \!\circ\! (\cD\tensor{B}\fk{a}) \!\!\!\!\!&\circ&\!\!\!\!\!(\Delta_{\cD}\tensor{B}\Sigma)
  =  (\cD\tensor{B}\tfks) \!\circ\! (\Delta_{\cD}\tensor{B}\tSi) \!\circ\!
  (\cD\tensor{B}\fk{a}) \!\circ\! (\Delta_{\cD}\tensor{B}\Sigma) \\
&=&\!\!\! (\cD\tensor{B}(\, \tfks \!\circ\! (\cD\tensor{B}\fk{a}) \!\circ\!(\Delta_{\cD}\tensor{B}\Sigma)\,))
     \!\circ\! (\Delta_{\cD}\tensor{B}\Sigma)  \\
   &=& \!\!\! (\cD\tensor{B}\fk{a}\tensor{A}\cC) \!\circ\!
   (\cD\tensor{B}\cD\tensor{B}\fk{s})
  \!\circ\! (\cD\tensor{B}\Delta_{\cD}\tensor{B}\Sigma) \!\circ\! (\Delta_{\cD}\tensor{B}\Sigma)
   \\
   &=&\!\!\! (\cD\tensor{B}\fk{a}\tensor{A}\cC) \!\circ\!
   (\cD\tensor{B}\cD\tensor{B}\fk{s})
  \!\circ\! (\Delta_{\cD}\tensor{B}\cD\tensor{B}\Sigma) \!\circ\! (\Delta_{\cD}\tensor{B}\Sigma) \\
   &=&\!\!\! (\cD\tensor{B}\fk{a}\tensor{A}\cC) \!\circ\!
   (\Delta_{\cD}\tensor{B}\Sigma\tensor{A}\cC) \!\circ\!
   (\cD\tensor{B}\fk{s})  \!\circ\! (\Delta_{\cD}\tensor{B}\Sigma) \\
   &=&\!\!\! ((\,(\cD\tensor{B}\fk{a}) \!\circ\!
   (\Delta_{\cD}\tensor{B}\Sigma)\,)\tensor{A}\cC) \!\circ\!
   \varrho^{\cD\tensor{B}\Sigma},
\end{eqnarray*}
thus proving the right $\cC$-colinearity of
$((\cD\tensor{B}\fk{a}) \circ
   (\Delta_{\cD}\tensor{B}\Sigma))$.

Conversely, consider a
$(\cD,\cC)$-bicolinear map $f: \cD\tensor{B}\Sigma \to \cD\tensor{B}\tSi$. We need to check that the bilinear map $\fk{a} (\varepsilon_{\cD}\tensor{B}\tSi) \circ f: \cD\tensor{B}\Sigma \to
\tSi$ is a $2$-cell in $\REM$. Using the $\cD$-colinearity of $f$, we compute
\begin{eqnarray*}
  \tfks \circ (\cD\tensor{B}\fk{a}) \circ (\Delta_{\cD}\!\tensor{B}\Sigma)
  &=& \tfks \circ (\cD\tensor{B}\varepsilon_{\cD}\tensor{B}\tSi)
  \circ (\cD\tensor{B}f) \circ (\Delta_{\cD}\tensor{B}\Sigma) \\
   &=& \tfks \circ (\cD\tensor{B}\varepsilon_{\cD}\tensor{B}\tSi)
  \circ (\Delta_{\cD}\tensor{B}\tSi) \circ f
   = \tfks \circ f
\end{eqnarray*}
On the other hand, \begin{multline*}
  (\fk{a}\tensor{A}\!\cC) \!\circ\! (\cD\tensor{B}\!\fk{s}) \!\circ\!
  (\Delta_{\cD}\tensor{B}\!\Sigma)
  = (\varepsilon_{\cD}\tensor{B}\tSi\tensor{A}\cC) \!\circ\! (f\tensor{A}\cC)
  \!\circ\! (\cD\tensor{B}\fk{s}) \!\circ\! (\Delta_{\cD}\tensor{B}\Sigma) \\
   \,\,=\,\, (\varepsilon_{\cD}\tensor{B}\tSi\tensor{A}\cC)
   \!\circ\! (\cD\tensor{B}\tfks) \!\circ\! (\Delta_{\cD}\tensor{B}\tSi) \!\circ\! f
   \,\,=\,\, \tfks \!\circ\! f,
\end{multline*}
where we have used the $\cC$-colinearity of $f$. Therefore,
$$\tfks \circ (\cD\tensor{B}\fk{a}) \circ
(\Delta_{\cD}\tensor{B}\Sigma) \,\, = \,\, (\fk{a}\tensor{A}\cC)
\circ (\cD\tensor{B}\fk{s}) \circ
(\Delta_{\cD}\tensor{B}\Sigma),$$
i.e.\ $\fk{a}$ is a $2$-cell as required.
\end{proof}

\begin{remark}
Morphisms between corings over different base rings have a natural meaning in $\mathsf{REM(Bim)}$. Given an algebra morphism $\alpha:B \to A$, any
$1$-cell from $(\cC\! :\! A)$ to $(\cD\! :\! B)$ of the form
$(A,\fk{s})$,  defines a morphism of corings $(\gamma,\alpha): (\cD\! :\! B)
\to (\cC\! :\! A)$ with
$$\xymatrix@C=50pt{\gamma: \cD \ar@{->}^-{\cD\tensor{B}\alpha}[r] & \cD\tensor{B}A
\ar@{->}^-{\fk{s}}[r] & A\tensor{A}\cC \simeq \cC.}
$$
 Conversely, any coring morphism
$(\gamma,\alpha):(\cD\! :\! B) \to (\cC\! :\! A)$ entails a $1$-cell
$(A,\fk{s})$ from $(\cC\! :\! A)$ to $(\cD\! :\! B)$ with
$$\xymatrix@C=50pt{\fk{s}: \cD\tensor{B}A
\ar@{->}^-{\gamma\tensor{B}A}[r] & \cC\tensor{B}A
\ar@{->}^-{\iota}[r] & \cC \simeq A\ot_A\cC,}$$
 where $\iota$
is the multiplication map.
\end{remark}

\subsection{The locally finite duality}
To study functors between categories of
comodules, it is convenient to introduce a different bicategory of
corings. Thus we define $\fREM$ as a bicategory obtained from
$\REM$ by restricting the class of $1$-cells to those that are
finitely generated and projective as right modules. Explicitly,
$\fREM$ has the same objects and 2-cells as $\REM$, while a 1-cell
$(\Sigma,\fk{s})$ from $(\cC\! :\! A)$ to $(\cD\! :\! B)$ in
$\REM$ is a 1-cell in $\fREM$ provided $\Sigma_A$ is a finitely
generated and projective module. The hom-category consisting of
$1$-cells from $(\cC\! :\! A)$ to $(\cD\! :\! B)$ and their
$2$-cells in the bicategory $\fREM$ is  denoted by ${}_{(\cD : B)}
{\mathcal{R}^{\mathsf{f}}}_{(\cC : A)}$.

In order to understand better the meaning of $\fREM$ for $\REM$ we
need to study the {\em left} bicategory of corings
$\LEM\,:=\, \mathsf{EM}\left(
\Bim^{co}\right)^{coop}$ (cf.\   \cite[p.\ 249]{Lack/Street:2002}).
$\LEM$ has corings as objects. $1$-cells from $(\cC\! :\! A)$ to
$(\cD\! :\! B)$ are pairs $(\Xi,\fk{x})$ consisting of an
$(A,B)$-bimodule $\Xi$ and an $(A,B)$-bilinear map $\fk{x}:
\Xi\ot_B\cD \rightarrow \cC\ot_A\Xi$, which is compatible with the
comultiplications and counits of both $\cD$ and $\cC$. The
identity $1$-cell associated to $(\cC\! :\! A)$ is the pair
$(A,\cC)$.
 If $(\Xi,\fk{x})$ and
$(\widetilde{\Xi},\widetilde{\fk{x}})$ are $1$-cells from $(\cC\!
:\! A)$ to $(\cD\! :\! B)$, then a $2$-cell is an $(A,B)$-bilinear
map $\fk{a}: \Xi\ot_B\cD \rightarrow \widetilde{\Xi}$, which is
compatible with $\fk{x}$, $\widetilde{\fk{x}}$ and the
comultiplication of $\cD$.

The category consisting of all $1$ and $2$-cells form $(\cC\! :\! A)$
to $(\cD\! :\! B)$ in $\LEM$ is denoted by
${}_{(\cD : B)}\mathcal{L}_{(\cC : A)}$. Furthermore,
${}_{(\cD : B)}{\mathcal{L}^{\mathsf{f}}}_{(\cC : A)}$  denotes the
full subcategory of ${}_{(\cD : B)}\mathcal{L}_{(\cC : A)}$ whose
objects are $1$-cells $(\Xi,\fk{x})$ such that ${}_A\Xi$ is
a finitely generated and projective left module. We use the notation
$\fLEM$ for the bicategory induced by the hom-categories
${}_{(\cD : B)}{\mathcal{L}^{\mathsf{f}}}_{(\cC : A)}$.

Recall that a \emph{duality} between categories is an
equivalence of categories via contravariant functors. The standard
duality between left and right  finitely
generated and projective modules induces a duality between the
hom-categories of $\fLEM$ and $\fREM$.

\begin{lemma}\label{duality}
Let $(\cC\! :\! A)$ and $(\cD\! :\! B)$ be corings. For any object $(\Sigma,\fk{s})$ in ${}_{(\cD : B)}{\mathcal{R}^{\mathsf{f}}}_{(\cC : A)}$, define $\fk{s}_*: \Sigma^*\ot_B\cD \to \cC\ot_A\Sigma^*$ by
$$
s^* \ot_B d \longmapsto \sum_i \left(\, (s^*\ot_A\cC) \circ
\fk{s}(d\ot_Be_i) \,\right)\ot_A e_i^*.
$$
For any morphism $\fk{a}:\cD\ot_B\Sigma \to \tSi$ in ${}_{(\cD : B)}{\mathcal{R}^{\mathsf{f}}}_{(\cC : A)}$, define $\fk{a}_*:
(\tSi)^*\ot_B\cD \rightarrow \Sigma^*$  by
$$\widetilde{s}^*\ot_Bd \longmapsto \sum_i
\widetilde{s}^*(\fk{a}(d\ot_Be_i))e_i^*.$$
The functor $$\xymatrix{ \mathscr{D}\left(\,(\cD\! :\! B),(\cC\! :\! A) \,\right):\,\,
{}_{(\cD : B)}{\mathcal{R}^{\mathsf{f}}}_{(\cC : A)} \ar@{->}[rr] & &
{}_{(\cD : B)}{\mathcal{L}^{\mathsf{f}}}_{(\cC : A)} },$$
given by  $\left( \, (\Sigma,\fk{s}) \mapsto (\Sigma^*,\fk{s}_*)\,
\right)$ and $\left(\, \fk{a} \mapsto \fk{a}_* \,\right)$ is
 a duality  of categories.
\end{lemma}
\begin{proof}
The maps $\fk{s}_*$, $\fk{a}_*$ are well-defined because, for all $b\in B$,
$\sum_i be_i\ot_A e^*_i = \sum_i e_i\ot_A e^*_ib$ and the
canonical element $\sum_i e_i\ot_A e^*_i$ is basis-independent.
The quasi-inverse contravariant functor is analogously constructed (its
action is denoted by the asterisk on the left). For every
object $(\Sigma, \fk{s}) \in
{}_{(\cD : B)}{\mathcal{R}^{\mathsf{f}}}_{(\cC : A)}$, we need to show
that the evaluating isomorphism $ev:\Sigma \simeq {}^*(\Sigma^*)$
is  an isomorphism in
${}_{(\cD : B)}{\mathcal{R}^{\mathsf{f}}}_{(\cC : A)}$. This follows from the
(easily checked) commutativity of the following diagram
$$\xymatrix@C=60pt@R=30pt{ \cD\ot_B \Sigma \ar@{->}^-{\fk{s}}[r]
\ar@{->}_-{\cD\ot_Bev}[d] &
\Sigma\ot_A\cC \ar@{->}^-{ev\ot_A\cC}[d] \\
\cD\tensor{B}{}^*(\Sigma^*) \ar@{->}^-{{}_*(\fk{s}_*)}[r] &
{}^*(\Sigma^*)\tensor{A}\cC. }$$
Therefore, $(\Sigma,\fk{s})
\simeq ({}^*(\Sigma^*), {}_*(\fk{s}_*))$  in
${}_{(\cD : B)}{\mathcal{R}^{\mathsf{f}}}_{(\cC : A)}$
by Proposition~\ref{red-unred}. The proof is completed
by routine computations that use the dual basis
criterion.
\end{proof}

Observe that the compatibility of the functors $\mathscr{D}(-,-)$
of Lemma \ref{duality} with the horizontal and vertical
compositions is guaranteed by the functoriality of tensor
products. In other words, functors $\mathscr{D}(-,-)$ give local
equivalences between the bicategories $\fREM$ and $(\fLEM)^{co}$.
Since both bicategories have the same objects, we obtain the
following

\begin{proposition}\label{Inv-morphBicat}
 The functors
$\mathscr{D}(-,-)$ of Lemma \ref{duality} establish a
biequivalence of bicategories $\fREM$ and $(\fLEM)^{co}$.
\end{proposition}

\section{Base ring extensions by modules}
\label{sec.base}

Base ring extensions of corings described in the introduction
correspond to extensions of base rings, i.e.\ to ring maps. Since
the work of Sugano \cite{Sug:not}, it has become clear that a more
general and unifying framework for studying ring extensions is
provided by bimodules rather than ring maps. In this section we
describe base extension of corings provided by a bimodule and also
study properties of modules reflected by base extensions of
corings.

\subsection{Definition of a base ring extension by a bimodule}
The basic idea of the construction of a base ring extension of corings by
a bimodule hinges on the relationship between comatrix corings and Sweedler's corings.

\begin{theorem}\label{thm.exten}
Given a coring $(\cD\! : \!B)$, let $\Sigma$ be a $(B,A)$-bimodule that
is finitely generated and projective as a right $A$-module.
Then the $A$-bimodule
$$
\Sigma[\cD] := \Sigma^*\ot_B\cD\ot_B \Sigma
$$
is an $A$-coring with the comultiplication
$$
\DSD : \Sigma[\cD]\to\Sigma[\cD]\ot_A\Sigma[\cD], \quad s^*\ot_Bd\ot_B s
\mapsto \sum_i s^*\ot_Bd\sw 1\ot_Be_i\ot_A e^*_i\ot_B d\sw 2\ot_B s,
$$
and the counit
$$
\eSD : \Sigma[\cD]\to A, \qquad s^*\ot_Bd\ot_B s \mapsto s^*(\eD(d)s).
$$
The coring $\Sigma[\cD]$ is called a {\em  base ring extension of a coring  by a module}.
\end{theorem}
\begin{proof}
A $B$-coring $\cD$ induces a comonad
$(F=-\tensor{B}\cD,\delta=-\tensor{B}\Delta_{\cD},\xi=-\tensor{B}\varepsilon_{\cD})$
in  $\rmod{B}$. On the other hand, $\Sigma$ induces an
adjunction
$$\xymatrix{S=-\tensor{A}\Sigma: \rmod{B}
\ar@<0,5ex>[r] & \ar@<0,5ex>[l] \rmod{A}:T=-\tensor{A}\Sigma^*}$$
with unit $\zeta$ and counit $\chi$,
given by, for all $X\in \cM_{B}$, $Y\in \cM_{A}$,
\begin{equation}\label{adj-uni-couni}
 \xymatrix@R=0pt{\zeta_{X}: X \ar@{->}[r] &
X\ot_B\Sigma\ot_A\Sigma^*, \\ x \ar@{|->}[r] & \sum_i
x\ot_Be_i\ot_Ae_i^*,} \quad \xymatrix@R=0pt{ \chi_{Y}:
Y\ot_A\Sigma^*\ot_B\Sigma \ar@{->}[r] & Y, \\ y\ot_As^*\ot_B s
\ar@{|->}[r] & ys^*(s). }
\end{equation}
 By \cite[Theorem 4.2]{Huber:1961}
(cf.\ the dual version of
\cite[Proposition 2.3]{Eilenberg/Moore:1965}),  these data give rise
to a comonad in
 $\rmod{A}$,
 $$(SFT, S F \zeta_{FT} \circ S\delta_{T}, \chi \circ S \xi_T).$$
Therefore,
$SFT(A)=A\tensor{A}\Sigma^*\tensor{B}\cD\tensor{B}\Sigma \simeq \Sigma^*\tensor{B}\cD\tensor{B}\Sigma$ is an $A$-coring.
The resulting comultiplication and counit
come out as stated.
\end{proof}

Examples of base ring extensions of a coring  by a module include
both base ring extensions and comatrix corings. In the former
case, given an algebra map $\alpha:B\to A$, one views $A$ as a
$(B,A)$-bimodule via $\alpha$ and multiplication in $A$, i.e.\
$baa' = \alpha(b)aa'$, for all $a,a'\in A$, $b\in B$. Obviously,
$A$ is a finitely generated and projective right $A$-module, and
the identification of $A$ with its dual $A^*$ immediately shows
that if $\Sigma = A$ then $\Sigma[\cD] \simeq A_\alpha[\cD]$. To
obtain a comatrix coring $\Sigma^*\ot_B\Sigma$, take $\cD$ to be
the trivial $B$-coring $B$.

\begin{remark}\label{coendomorph}
By \cite[Example 5.1]{KaoGom:com}, comatrix
corings can be understood as the coendomorphism corings of suitable
quasi-finite bicomodules. The coring introduced in Theorem~\ref{thm.exten} can
also be understood in this way. Start
with the following adjunctions
$$
\xymatrix@R=30pt@C=60pt{ \rcomod{\cD} \ar@<0,5ex>^-{U_B}[r]
& \rmod{B} \ar@<0,5ex>^-{-\tensor{B}\cD}[l]
\ar@<0,5ex>^-{-\tensor{B}\Sigma}[rd] & \\
& & \ar@<0,5ex>^-{-\tensor{A}\Sigma^*}[lu] \rmod{A}, }
$$
where $U_B$ is the forgetful functor. The counit and unit of the
composite adjunction
\begin{equation}\label{eq.adj1}
\xymatrix@C=60pt{ -\tensor{A}(\Sigma^*\tensor{B}\cD): \rmod{A}
\ar@<0,5ex>[r] &  \rcomod{\cD}: (-\tensor{B}\Sigma) \circ U_B
\ar@<0,5ex>[l],  }
\end{equation}
are given, respectively, for all $Y\in \rmod{A}$ and $X\in \rcomod{\cD}$,
\begin{equation}\label{eq.counit}
\xymatrix@R=0pt{ \psi_{Y}:
U_B(Y\tensor{A}\Sigma^*\tensor{B}\cD)\tensor{B}\Sigma \ar@{->}[r]
& Y, \\
\sum_{\alpha}y^{\alpha}\tensor{A}s_{\alpha}^*\tensor{B}d^{\alpha}\tensor{B}u_{\alpha}
\ar@{|->}[r] & \sum_{\alpha}y^{\alpha}s^*_{\alpha}\left(
\eps_{\cD}(d^{\alpha})u_{\alpha}\right) }
\]
and
\begin{equation}\label{eq.unit}
\xymatrix@R=0pt{ \eta_{X}: X \ar@{->}[r] &
\left(U_B(X)\tensor{B}\Sigma\right)\tensor{A}\Sigma^*\tensor{B} \cD ,\\
x \ar@{|->}[r] & \sum_{i,(x)}x_{(0)}\tensor{B} e_i \tensor{A}
e_i^* \tensor{B} x_{(1)}. }
\end{equation}
Since $(-\tensor{B}\Sigma) \circ U_B$ is the left
adjoint to $-\tensor{A}(\Sigma^*\tensor{B}\cD)$,
$\Sigma^*\tensor{B}\cD$ is an $(A,\cD)$-quasi-finite
$(A,\cD)$-bicomodule, where $A$ is considered trivially as an
$A$-coring (cf.\ \cite[Definition 4.1]{Gom:sep}). The corresponding co-hom
functor comes out as $ h_{\cD}(\Sigma^*\tensor{B}\cD, -) U_B(-)\tensor{B}\Sigma$. Therefore, the discussion of
\cite[Section 5]{KaoGom:com} implies that the coendomorphism
$A$-coring associated to
$\Sigma^*\tensor{B}\cD$ is
\[
e_{\cD}(\Sigma^*\tensor{B}\cD) h_{\cD}(\Sigma^*\tensor{B}\cD,
\Sigma^*\tensor{B}\cD) \Sigma^*\tensor{B}\cD\tensor{B}\Sigma.
\]
This is precisely the coring constructed in Theorem
\ref{thm.exten}, and the forms of the unit and the counit of the
defining adjunction imply immediately the stated form of the
comultiplication and counit.
\end{remark}

The constructions given in Theorem \ref{thm.exten} are
functorial. The following proposition summarises  basic
properties of these functors.

\begin{proposition}\label{pro.functor}
\label{functor}
Let $\Sigma$ be a $(B,A)$-bimodule that is finitely generated
and projective as a right $A$-module.
\begin{zlist}
\item  The assignment
$$
\cD\mapsto \Sigma[\cD], \qquad f\mapsto \Sigma^*\ot_B f\ot_B \Sigma
$$
defines a functor $\Sigma[-]$ : $B$-\cor $\to$ $A$-\cor ,\; which
commutes with colimits. \item For any $(C,B)$-bimodule $\Xi$ that
is finitely generated and projective as a right $B$-module,
$$
\Sigma[-]\circ \Xi [-] \simeq (\Xi\ot_B\Sigma)[-]
$$
as functors.
\end{zlist}
\end{proposition}
\begin{proof}
(1) A straightforward calculation that uses the definitions of the
comultiplication and counit in $\Sigma[\cD]$ and the fact that $f$
is a morphism of $B$-corings, confirms that $\Sigma^*\ot_B f\ot_B
\Sigma$ is a morphism of $A$-corings. Since the compositions in
$B$-\cor{} and $A$-\cor{} are the same as composition of
$k$-modules, $\Sigma[-]$ is a functor. That $\Sigma[-]$ commutes
with colimits is an easy consequence of the fact that tensor
products commute with colimits, and that the colimit of an
inductive system of $B$-corings is already computed in the
category of $B$-bimodules.

(2) For any $C$-coring $\cD$, consider the $(A,A)$-bimodule isomorphism
\[
\phi_{\Sigma,\Xi,\cD}: \Sigma ^* \ot_B\Xi^*\ot_C\cD\ot_C\Xi\ot_B \Sigma
\to  (\Xi\tensor{B}\Sigma)^*\ot_C\cD\ot_C(\Xi\tensor{B}\Sigma),\]
\[  s^*\ot_B x^*\ot_C d\ot_C x\ot_B s \mapsto \overline{(s^*\ot_B x^*)}\ot_C d\ot_C x\ot_B s,
\]
where $\bara{(-)}: \Sigma^* \tensor{B} \Xi^*  \rightarrow
(\Xi\tensor{B}\Sigma)^*$ is the $(A,C)$-bilinear isomorphism sending $s^*\ot_B x^*$ to
$x \tensor{B}s \mapsto \bara{(s^*\tensor{B}x^*)}(x
\tensor{B}s) =s^*(x^*(x)s)$.
We need to prove that $\phi_{\Sigma,\Xi,\cD}$ is an
$A$-coring map, natural in $\cD$. Write $V=\Xi\tensor{B}\Sigma$.  Clearly $V$ is a
$(C,A)$-bimodule that is  finitely generated and
projective as a right $A$-module. Let
$\{f_k,f_k^*\}$ be a dual basis for $\Xi_B$ and let $\{e_i,e_i^*\} $ be a dual basis for $\Sigma_A$.
Then $\{f_k\tensor{B}e_i, \bara{(e_i^*\tensor{B}f_k^*)} \}_{i,k}$
is a finite  dual basis for $V_A$. For all  $s^*\in
\Sigma^*$, $x^* \in \Xi^*$, $x \in \Xi$, $s \in \Sigma$ and $d\in \cD$,
\begin{eqnarray*}
  \eps_{V[\cD]}(\phi_{\Sigma,\Xi,\cD}(s^*\tensor{B}x^*\tensor{C}d\tensor{C}
  x\tensor{B}s))
   &=& \bara{(s^*\tensor{B} x^*)}(\eD(d)x\tensor{B}s)
   \\   &=& s^*(x^*(\eD(d)x)s) \\
   &=& \eps_{\Sigma[\Xi[\cD]]}(s^*\tensor{B}x^*\tensor{C}d\tensor{C}
  x\tensor{B}s).
\end{eqnarray*}
Furthermore,
\begin{eqnarray*}
  \Delta_{V[\cD]} \!\!\!\!\!\!\!\!\!\!\!\!\!\!\!&&(\phi_{\Sigma,\Xi,\cD}
  (s^*\tensor{B}x^*\tensor{C}d\tensor{C}
  x\tensor{B}s)) \\
   &=&\sum_{i,k}\bara{(s^*\tensor{B} x^*)} \tensor{C}d\sw 1\tensor{C}
   (f_k\tensor{B}e_i) \tensor{A} \bara{(e_i^*\tensor{B}f_k^*)} \tensor{C}
   d\sw 2\tensor{C}(x\tensor{B}s) \\
   &=& \left( \phi_{\Sigma,\Xi,\cD} \!\tensor{A}\! \phi_{\Sigma,\Xi,\cD} \right) \!
   \left(\sum_{i,k} s^*\!\tensor{B}\! x^* \!\tensor{C}\! d\sw 1\!\tensor{C}\!
   f_k\!\tensor{B}\!e_i \!\tensor{A}\! e_i^*\!\tensor{B}\! f_k^*
   \!\tensor{C}\! d\sw 2\! \tensor{C}\! x\!\tensor{B}\! s \right) \\
   &=& \left( \phi_{\Sigma,\Xi,\cD} \tensor{A} \phi_{\Sigma,\Xi,\cD} \right)
    \circ \Delta_{\Sigma[\Xi[\cD]]}
     (s^*\tensor{B}x^*\tensor{C}d\tensor{C}
  x\tensor{B}s).
\end{eqnarray*}
This proves that
\[
\eps_{(\Xi\ot_B\Sigma)[\cD]}\circ \phi_{\Sigma,\Xi,\cD} \eps_{\Sigma[\Xi[\cD]]}, \text{ and }
\Delta_{(\Xi\ot_B\Sigma)[\cD]}\circ \phi_{\Sigma,\Xi,\cD} = \left(
\phi_{\Sigma,\Xi,\cD} \tensor{A} \phi_{\Sigma,\Xi,\cD} \right)
    \circ \Delta_{\Sigma[\Xi[\cD]]},
\]
i.e.\ that for all $C$-corings $\cD$, the $(A,A)$-bimodule
isomorphism   $\phi_{\Sigma,\Xi,\cD}$ is a morphism of $A$-corings.
Since a bijective morphism of corings is an isomorphism,
all the $\phi_{\Sigma,\Xi,\cD}$ are isomorphisms of $A$-corings
and their explicit forms immediately imply that they are natural in $\cD$.
\end{proof}

In particular, Proposition~\ref{functor} leads to
\begin{corollary}\label{eps.tilde}
Let $\Sigma$ be a $(B,A)$-bimodule  that is finitely generated and
projective as a right $A$-module
and let $\cD$ be a $B$-coring.
\begin{enumerate}
\item The map
$$\eps_{\cD,\Sigma}: \Sigma[\cD] \to \Sigma^*\otimes_B \Sigma,
\quad s^*\ot_B d\ot_B s\mapsto s^*\ot_B\eD(d)s,
$$
is a morphism of $A$-corings.

\item The maps $$\fk{s}: \cD\ot_B\Sigma \to \Sigma \ot_A
\Sigma^*\ot_B\Sigma, \quad d \ot_B s \mapsto \sum_ie_i\ot_Ae_i^*
\ot_B \eps_{\cD}(d)s,$$ and $$\fk{s}':\cD\ot_B\Sigma \to
\Sigma\ot_A\Sigma[\cD],\quad d \ot_B s \mapsto \sum_ie_i\ot_Ae_i^*
\ot_B d\ot_Bs,$$ define, respectively, $1$-cells $(\Sigma,\fk{s})$
and $(\Sigma,\fk{s}')$ from $(\Sigma^*\ot_B \Sigma\! :\! A)$ to
$(\cD\! :\! B)$ and from $(\Sigma[\cD]\! :\! A)$ to $(\cD\! :\! B)$ in the
bicategory $\fREM$.

\item The maps $$\fk{t}:\Sigma^*\ot_B \cD \to \Sigma^*\ot_B\Sigma
\ot_A \Sigma^*,\quad s^*\ot_B d \mapsto \sum_is^*\eps_{\cD}(d)
\ot_Be_i\ot_Ae_i^*,$$ and $$\fk{t}': \Sigma^*\ot_B\cD \to
\Sigma[\cD]\ot_A\Sigma^*,\quad s^*\ot_Bd \mapsto
\sum_is^*\ot_Bd\ot_Be_i\ot_Ae_i^*,$$ define, respectively,
$1$-cells $(\Sigma^*,\fk{t})$ and $(\Sigma^*,\fk{t}')$ from
$(\Sigma^*\ot_B \Sigma\! :\! A)$ to $(\cD\! :\! B)$ and from $(\Sigma[\cD]\! :\!
A)$ to $(\cD\! :\! B)$ in the bicategory $\fLEM$.
\end{enumerate}
\end{corollary}
\begin{proof}
$(1)$ View $B$ as a trivial $B$-coring. Then $\eD:\cD\to B$ is a
morphism of $B$-corings. Note that $\eps_{\cD,\Sigma} \Sigma[\eD]$, hence it is a morphism of $A$-corings. The proofs
of statements $(2)$ and $(3)$ use the dual basis criterion, and are left
 to the reader.
\end{proof}

Given a comatrix coring $\Sigma^*\ot_B\Sigma$, $\Sigma$ is a right comodule and $\Sigma^*$ is a left comodule of $\Sigma^*\ot_B\Sigma$
 (cf.\ \cite[p.\ 891]{KaoGom:com}).
 This can be extended to the following observation.

\begin{lemma}
\label{comodule}
Given a coring $(\cD\! : \!B)$, let $\Sigma$ be a $(B,A)$-bimodule
that is finitely generated and projective as a right $A$-module and let
$\Sigma[\cD]$ be the associated base ring extension coring.
Then:
\begin{zlist}
\item  $\cD\ot_B\Sigma$ is a  $(\cD,\Sigma[\cD])$-bicomodule with
the left coaction ${}^{\cD\ot_B\Sigma}\varrho = \DD\ot_B \Sigma$
and the right coaction
$$
\varrho^{\cD\ot_B\Sigma}: \cD\ot_B\Sigma\to \cD\ot_B\Sigma\ot_A \Sigma[\cD],
\quad d\ot_B s\mapsto \sum_i d\sw 1\ot_B e_i\ot_A e^*_i\ot_B d\sw 2\ot_B s.
$$

\item $\Sigma^*\ot_B\cD$ is a $(\Sigma[\cD],\cD)$-bicomodule
with the right coaction $\varrho^{\Sigma^*\ot_B\cD}= \Sigma^*\ot_B \DD$
and the left coaction
$$
{}^{\Sigma^*\ot_B\cD}\varrho: \Sigma^*\ot_B\cD \to \Sigma[\cD]\ot_A\Sigma^*\ot_B \cD,
\quad s^*\ot_B d\mapsto \sum_i s^*\ot_B d\sw 1\ot_B e_i\ot_A e^*_i\ot_B d\sw 2.
$$
\end{zlist}
\end{lemma}
\begin{proof}
This is a direct consequence of Proposition~\ref{push-Out} (and its left-handed
version) and Corollary~\ref{eps.tilde}.
\end{proof}

\subsection{Properties reflected by base ring extensions}
In \cite{BrzGom:com} it has been studied how properties
of bimodules are reflected by the  associated finite comatrix
corings. The aim of this section is to study such reflection
of properties by the base ring
extensions of corings by modules and thus
to generalise the main results of \cite{BrzGom:com}.

Following \cite{Sug:not}, a $(B,A)$-bimodule $\Sigma$ is said
to be a \emph{separable bimodule}, if the evaluation map
\[
\Sigma\tensor{A}{}^*\Sigma \rightarrow B, \quad s\tensor{A}s^*
\mapsto s^*(s)
\]
is a split epimorphism of $(B,B)$-bimodules. Furthermore, recall from
\cite{AndFul:rin} and \cite{Kad:new}  that $\Sigma$ is
said to be a \emph{Frobenius bimodule} if ${}_B\Sigma$ and
$\Sigma_A$ are finitely generated projective modules, and
$\Sigma^* \simeq {}^*\Sigma$ as $(A,B)$-bimodules.

An $A$-coring
$\coring{C}$ is said to be a \emph{Frobenius coring} if the
forgetful functor $U_A: \rcomod{\coring{C}} \rightarrow \rmod{A}$
is a Frobenius functor in the sense of
\cite{CaeMil:Doi} and
\cite{CasGom:Fro} (cf.\
\cite{Brz:str} and \cite{Brz:tow} for more
details and other equivalent characterisations of a Frobenius coring).
An $A$-coring $\coring{C}$ is
said to be  a \emph{cosplit coring} if and only if
$\eps_{\coring{C}}$ is a split epimorphism of $(A,A)$-bimodules.
Equivalently, $\cC$ is a cosplit coring if there exists an invariant element $c \in
\coring{C}$ (i.e.\ $ac=ca$, for all $a \in A$) such that
$\eps_{\coring{C}}(c)=1_A$. An $A$-coring $\coring{C}$ is called  a
\emph{coseparable coring} if and only if $\Delta_{\coring{C}}$ is
a split monomorphism of $\coring{C}$-bicomodules (cf.\ \cite{Guz:coi} and
\cite{GomLou:cos} for more details on coseparable corings).

Take a $(B,A)$-bimodule $\Sigma$ that is
 finitely generated and projective as a right $A$-module,
and consider its right endomorphism ring
$S=\rend{A}{\Sigma}$. Then there is a canonical isomorphism of
$(S,S)$-bimodules $\Sigma\tensor{A}\Sigma^* \simeq S$ given by
 $\Sigma\tensor{A}\Sigma^*\ni s\tensor{A}s^* \to [s' \mapsto ss^*(s')]$.
 Its inverse is given by $f\mapsto \sum_i
f(e_i)\tensor{A}e_i^*$.

\begin{proposition}
Let $\coring{D}$ be a $B$-coring and $\Sigma$ a $(B,A)$-bimodule
such that $\Sigma_A$ is a finitely generated projective module.
\begin{zlist}
\item If ${}_B\Sigma_A$ is a separable bimodule and $\coring{D}$
is a coseparable $B$-coring, then $\Sigma[\coring{D}]$ is a
coseparable $A$-coring.

\item If ${}_A\Sigma^*_B$ is a separable bimodule and
$\coring{D}$ is a cosplit $B$-coring, then $\Sigma[\coring{D}]$
is a cosplit $A$-coring. Conversely, if $\Sigma[\coring{D}]$
is a cosplit $A$-coring, then ${}_A\Sigma^*_B$ is a separable
bimodule.

\item If ${}_B\Sigma_A$ is a Frobenius bimodule and $\coring{D}$
is a Frobenius $B$-coring, then $\Sigma[\coring{D}]$ is a
Frobenius $A$-coring.
\end{zlist}
\end{proposition}
\begin{proof}
(1) By \cite[Theorem 3.1(1)]{Kad:sep}, the ring extension $B
\rightarrow S$ is split. Let $\kappa: S \rightarrow B$ be a
$B$-bilinear splitting map. Since $1_S$ is sent to
$\sum_ie_i\tensor{A}e_i^*$ by the isomorphism $S \simeq
\Sigma\tensor{A}\Sigma^*$, $\kappa(\sum_ie_i\tensor{A}e_i^*)=1_B$.
Let $\nabla_{\coring{D}}: \coring{D}\tensor{B}\coring{D}
\rightarrow \coring{D}$ be a $\coring{D}$-bicolinear map such that
 $\nabla_{\coring{D}} \circ \Delta_{\coring{D}} \coring{D}$. Consider an  $A$-bilinear map,
\[
\xymatrix@R=0pt{ \nabla_{\Sigma[\coring{D}]}:
\Sigma[\coring{D}] \tensor{A} \Sigma[\coring{D}]
\ar@{->}[r] & \Sigma[\coring{D}], \\
s^*\tensor{B}d\tensor{B}s \tensor{A}
\ts^*\tensor{B}d'\tensor{B}\ts \ar@{|->}[r] &
s^*\tensor{B}\nabla_{\coring{D}}\left(d\kappa(s\tensor{A}\ts^*)\tensor{B}d'
\right)\tensor{B}\ts .}
\]
Take any $s^*\in \Sigma^*$, $s\in \Sigma$ and $d\in \cD$ and compute
\begin{eqnarray*}
  \nabla_{\Sigma[\coring{D}]} \circ \Delta_{\Sigma[\coring{D}]} ( s^*\tensor{B}d\tensor{B}s)
  &=& \nabla_{\Sigma[\coring{D}]}
  \left( \sum_{i} s^*\tensor{B}d_{(1)}\tensor{B}e_i
  \tensor{A} e_i^*\tensor{B}d_{(2)}\tensor{B}s \right)
  \\
   &=& \sum_{i} s^*\tensor{B}\nabla_{\coring{D}}(d_{(1)}
   \kappa(e_i\tensor{A}e_i^*)\tensor{B}d_{(2)})\tensor{B}s \\
   &=&
   \sum s^*\tensor{B}\nabla_{\coring{D}}(d_{(1)}\tensor{B}d_{(2)})\tensor{B}s
   \\ &=& s^*\tensor{B}d\tensor{B}s.
\end{eqnarray*}
Thus,  $\nabla_{\Sigma[\coring{D}]} \circ \Delta_{\Sigma[\coring{D}]} \Sigma[\coring{D}]$, as required. We need to show that
$\nabla_{\Sigma[\coring{D}]}$ is a
$\Sigma[\coring{D}]$-bicolinear map.  Take any $s^*,\ts^* \in
\Sigma^*$, $d,d' \in \coring{D}$ and  $s,\ts\in \Sigma$, and compute
\begin{multline*}
\Delta_{\Sigma[\coring{D}]} \circ \nabla_{\Sigma[\coring{D}]}
(s^*\tensor{B}d\tensor{B}s \tensor{A}
\ts^*\tensor{B}d'\tensor{B}\ts) \\ = \sum_{i} s^* \tensor{B}
\nabla_{\coring{D}}(d\kappa(s\tensor{A}\ts^*)\tensor{B}d_{(1)}'
)\tensor{B}e_i\tensor{A}e_i^*\tensor{B}d_{(2)}'\tensor{B}\ts
\\ (\nabla_{\Sigma[\coring{D}]}\tensor{A}{\Sigma[\coring{D}]})(\sum_{i}
s^*\tensor{B}d\tensor{B}s\tensor{A}
\ts^*\tensor{B}d'_{(1)}\tensor{B}e_i \tensor{A}e_i^*
\tensor{B}d_{(2)}'\tensor{B}\ts)
\\  (\nabla_{\Sigma[\coring{D}]}\tensor{A}{\Sigma[\coring{D}]})\circ
({\Sigma[\coring{D}]}\tensor{A}\Delta_{\Sigma[\coring{D}]})(s^*\tensor{B}d\tensor{B}s
\tensor{A} \ts^*\tensor{B}d'\tensor{B}\ts).
\end{multline*}
To derive this equality we have used the right $\coring{D}$-colinearity of
$\nabla_{\coring{D}}$. This proves that  $\nabla_{\Sigma[\coring{D}]}$ is
a right ${\Sigma[\coring{D}]}$-colinear map. Similarly, one uses the
left $\coring{D}$-colinearity of $\nabla_{\coring{D}}$ to obtain
the left ${\Sigma[\coring{D}]}$-colinearity of
$\nabla_{\Sigma[\coring{D}]}$. Thus we conclude that ${\Sigma[\coring{D}]}$ is
a coseparable coring as claimed.

(2) Consider the composition of $A$-bilinear maps $\theta:A
\rightarrow \Sigma^*\tensor{B}{}^*(\Sigma^*) \simeq
\Sigma^*\tensor{B}\Sigma$, where the first map is a splitting of the
evaluation map that is provided by the separability of
${}_A\Sigma^*_B$. Put $\theta(1_A) = \sum_k s^*_k \tensor{B}s_k$.
Note that $\theta(1_A)$ is an invariant element of the
$A$-bimodule $\Sigma^*\tensor{B}\Sigma$ that in addition satisfies
$\sum_ks^*_k(s_k)=1_A$. By hypothesis, $\coring{D}$ is a cosplit
$B$-coring, hence there exists an invariant  element $d \in
\coring{D}$ such that $\eps_{\coring{D}}(d)=1_B$. Combining $d$
with $\theta(1_A)$ one obtains an invariant element
 $\sum_ks^*_k\tensor{B}d\tensor{B}s_k \in
\Sigma[\coring{D}]$, such that
$\eps_{\Sigma[\coring{D}]}(\sum_ks^*_k\tensor{B}d\tensor{B}s_k )=1_A$.

Conversely, if $\Sigma[\coring{D}]$ is a cosplit $A$-coring, and
$\sum_{l}s^*_l \tensor{B}d_l\ot_Bs_l$ is the corresponding invariant element,
then $\sum_{l}s^*_l \tensor{B}\eps_{\coring{D}}(d_l)s_l
\in \Sigma^*\tensor{B}\Sigma$ gives a section for the evaluation
$\Sigma^*\tensor{B}{}^*(\Sigma^*) \rightarrow A$.

(3) Suppose that $\Sigma$ is a Frobenius $(B,A)$-bimodule and let
$\gamma: \Sigma^* \rightarrow {}^*\Sigma$ be the corresponding
$(A,B)$-bilinear isomorphism. Let $T$ be the opposite ring to the
left dual ring ${}^*\cD$. By \cite[Theorem 4.1]{Brz:str}, there is
an isomorphism $\sigma: \coring{D} \to T$ of $(B,T)$-bimodules and
${}_B\coring{D}$ is a finitely generated and projective module.
Consider an $A$-bilinear map $\beta({}^*\Sigma\tensor{B}\sigma\tensor{B}\Sigma) \circ
(\gamma\tensor{B}\coring{D} \tensor{B} \Sigma)$. This leads to an
$(A,B)$-bimodule isomorphism
\[
\xymatrix{  \zeta: \Sigma^* \tensor{B} \coring{D} \tensor{B}
\Sigma \ar@{->}^-{\beta}[r] &
{}^*\Sigma \tensor{B} {}^*\coring{D}
\tensor{B} \Sigma \ar@{->}^-{\simeq}[r] &
\lhom{B}{\coring{D}\tensor{B}\Sigma}{\Sigma}}
\]
where the last isomorphism is a consequence of the fact that both
${}_B\coring{D}$ and ${}_{B}\Sigma$ are finitely generated and
projective modules. Note that ${}_A\Sigma[\coring{D}]$ is a
finitely generated and projective module. In view of \cite[Theorem
4.1]{Brz:str}, it suffices to show that $\Sigma[\coring{D}] \simeq
R$ as  $(A,R)$-bimodules, where $R$ is the opposite ring of the
left dual ring  ${}^*\Sigma[\coring{D}]$. Consider the following
chain of isomorphisms
\[
\xymatrix{ \Sigma^*\tensor{B} \coring{D} \tensor{B}\Sigma
\ar@{->}^-{\zeta}[r] & \lhom{B}{\coring{D}\tensor{B}\Sigma}{\Sigma}
\ar@{->}^-{\simeq}[r] &
\lhom{B}{\coring{D}\tensor{B}\Sigma}{\lhom A{\Sigma^*}A}&\\
&\ar@{->}^-{\simeq}[r] & \lhom A {\Sigma^*\ot_B\cD\ot_B\Sigma}A = R, }
\]
where the second isomorphism follows from the fact that $\Sigma_A$ is
a finitely generated and projective module.
One easily checks that all these isomorphisms are $(A,R)$-bimodule maps and thus their
composition provides one with  the required map.
\end{proof}

\begin{remark}
There are properties of a $B$-coring $\cD$ that are directly
reflected in  $A$-coring $\Sigma[\cD]$ without any assumption
on  bimodule $\Sigma$. For instance, if $\cD$ is a cosemisimple
$B$-coring \cite{KaoGomLob:sem}, then by \cite[Theorem
4.4]{KaoGom:com}, $\cD$ decomposes into a  direct sum of comatrix
corings, $\cD \simeq \oplus_{\Xi \in
\Lambda}\rcomatrix{D_{\Xi}}{\Xi}$, where each of the $\Xi$ is a
finitely generated and projective right $B$-module and each
$D_{\Xi}$ is a division subring of the endomorphism ring
$\rend{B}{\Xi}$. Thus, $\cD \simeq \oplus_{\Xi \in \Lambda}
\Xi[D_{\Xi}]$ as $B$-corings. Applying the functor $\Sigma[-]$ to
this isomorphism and using Proposition \ref{pro.functor}, one
obtains an isomorphism of $A$-corings
\begin{eqnarray*}
\Sigma[\coring{D}] &\simeq & \oplus_{\Xi \in\Lambda} \Sigma[-] \circ \Xi[D_{\Xi}]  \\
&\simeq & \oplus_{\Xi \in \Lambda} (\Xi\tensor{B}\Sigma)[D_{\Xi}].
\end{eqnarray*}
Therefore, \cite[Theorem
4.4]{KaoGom:com} implies that $\Sigma[\cD]$ is a cosemisimple
$A$-coring.
\end{remark}

\section{Module-morphisms and push-out and pull-back
functors}\label{sec.mod}
Following the general strategy of replacing algebra maps by
bimodules we introduce the notion of a
{\em module-morphism} of corings and study
properties of associated functors between  categories of comodules.
\subsection{The categories of module-morphisms and
module-representations}\label{mod-moph-rep}

\begin{definition}\label{def.mod-mor}
Let $(\cC\! : \!A)$ and $(\cD\! : \!B)$ be corings. A {\em
$((\cD\! : \!B), (\cC\! : \!A))$-module-morphism}  is a pair $\bSi
= (\Sigma, \sigma)$ where $\Sigma$ is a $(B,A)$-bimodule that is
finitely generated and projective as a right $A$-module and
$\sigma: \Sigma[\cD] \to \cC$ is a morphism of $A$-corings. A
$((\cD\! : \!B), (\cC\! : \!A))$-module-morphism $\bSi$ is often
denoted by ${}_{(\cD : B)}\bSi_{(\cC : A)}$.
\end{definition}

\begin{example}\label{ex.mod-mor} \rm ~\\
(1)  Any morphism of corings $(\gamma,\alpha): (\cD\! : \!B)\to (\cC\! : \!A)$
gives rise to a module-morphism $(A,\tilde{\gamma})$, where $A$ is  a left $B$-module via
$\alpha$ and $\tilde{\gamma}: A\ot_B\cD\ot_B A\to \cC$ is the induced map
$a\ot_B d\ot_B a'\mapsto a\gamma(d)a'$.

\noindent (2)  Let $\cC$ be an
$A$-coring and let $M$ be a right $\cC$-comodule that is finitely generated and
projective as a right $A$-module. Let $B$ be the endomorphism ring
$B = \Rend\cC M$ so that $M$ is a $(B,A)$-bimodule. Take any $B$-coring $\cD$. Then
$$
\can_{\cD,M} : M^*\ot_B\cD\ot_B M\to \cC,\qquad m^*\ot_B d\ot_B m\mapsto
\sum m^*(\eD(d)m\sw 0)m\sw 1,
$$
is an $A$-coring map. Thus $(M,\can_{\cD,M})$ is a
$((\cD\! : \!B), (\cC\! : \!A))$-module-morphism.
\end{example}
\begin{proof}
Example (1) follows immediately from the definition of a morphism
of corings. To check (2) simply note that
$\can_{\cD,M} = \can_M\circ \eps_{\cD,M}$, where $\eps_{\cD,M}$ is
$A$-coring morphism in Corollary~\ref{eps.tilde} and $\can_M := \can_{B,M}$ is
an $A$-coring morphism by \cite[18.26]{BrzWis:cor} and
\cite[Proposition 2.7]{KaoGom:com}. Thus
$\can_{\cD,M}$ is an $A$-coring morphism as required.
\end{proof}

Corollary~\ref{eps.tilde} leads to the following interpretation of
module-morphisms in terms of
 $1$-cells in the bicategory $\fREM$,
\begin{lemma}\label{lemma.mod-mor}
$((\cD\! : \!B), (\cC\! : \!A))$-module-morphisms are in bijective correspondence with
1-cells from $(\cC\! :\! A)$ to
$(\cD\! :\! B)$ in $\fREM$.
\end{lemma}
\begin{proof}
The correspondence follows from  the
natural isomorphism
\begin{equation}\label{eq.1-iso}
\hom {A,A} {\Sigma^*\otimes_B\cD\otimes_B\Sigma}{\cC} \simeq
\hom{B,A}{\cD\otimes_B\Sigma}{\Sigma\ot_A\cC}.
\end{equation}
Explicitly, if   $(\Sigma, \sigma)$ is a $((\cD\! : \!B), (\cC\! : \!A))$-module-morphism,
 then  $(\Sigma,\fk{s}_{\sigma})$, with
 \begin{equation}\label{fks-sigma}
\fk{s}_{\sigma}:\cD\ot_B\Sigma \to \Sigma\ot_A\cC,\quad d\ot_Bs
\mapsto \sum_ie_i\ot_A \sigma(e_i^*\ot_Bd\ot_Bs)
\end{equation}
is a 1-cell in  $\fREM$.
Conversely, given a 1-cell $(\Sigma,\fk{s})$ in $\fREM$, define $\sigma_\fk{s}$ as
the composition
$$
\xymatrix@C=60pt{\sigma_{\fk{s}}: \Sigma^*\ot_B\cD\ot_B\Sigma
\ar@{->}^-{\Sigma^*\ot_B\fk{s}}[r] & \Sigma^*\ot_B\Sigma\ot_A\cC
\ar@{->}^-{\eps_{\Sigma[B]}\ot_A\cC}[r] & A\ot_A \cC \simeq \cC. }
$$
Then $(\Sigma, \sigma_{\fk{s}})$ is a module-morphism. That this
correspondence is well-defined and bijective can be checked directly by
using the properties of dual bases.
\end{proof}

\begin{remark}\label{left-1-cell}
Put together, Lemma~\ref{lemma.mod-mor} and
Lemma~\ref{duality}  establish a bijective correspondence between module-morphisms and  $1$-cells in the bicategory
$\fLEM$. For a module-morphism
${}_{(\cD : B)}\bSi_{(\cC : A)}$, define the $(A,B)$-bilinear map
\begin{equation}\label{fkt}
\fk{t}_{\sigma}:\Sigma^*\ot_B\cD \to \cC\ot_A\Sigma^*,\quad
s^*\ot_Bd \mapsto \sum_i \sigma(s^*\ot_B d \ot_B e_i) \ot_A e_i^*.
\end{equation}
Then $(\Sigma^*,\fk{t}_{\sigma})$ is
a $1$-cell from $(\cC\! :\! A)$ to
$(\cD\! :\! B)$ in  $\fLEM$. In the converse direction, an object $(\Xi,\fk{x})\in {}_{(\cD : B)}{\mathcal{L}^{\mathsf{f}}}_{(\cC : A)}$
induces a module-morphism $({}^*\Xi,  \sigma_{\fk{x}})$, where
$$
\xymatrix@C=40pt{ \sigma_{\fk{x}}: ({}^*\Xi)^* \ot_B \cD\ot_B{}^*\Xi \simeq
\Xi\ot_B\cD\ot_B{}^*\Xi \ar@{->}^-{\fk{t}\ot_B{}^*\Xi}[r] &
\cC\ot_A\Xi\ot_B{}^*\Xi \ar@{->}^-{\cC\ot_A\eps_{{}^*\Xi[B]}}[r] & \cC\ot_A A
\simeq \cC.}
$$
\end{remark}

In view of Proposition~\ref{red-unred} and Lemma~\ref{lemma.mod-mor},
 for any module-morphism
${}_{(\cD : B)}\bSi_{(\cC : A)}$, the $(B,A)$-bimodule $\cD\ot_B
\Sigma$ is a $(\cD,\cC)$-bicomodule with the left $\cD$-coaction
$\DD\ot_B\Sigma$ and the right $\cC$-coaction $\varrho^{\bSi_\circ
(\cD)}$ given by Proposition~\ref{push-Out}. A {\em map of
$((\cD\! : \!B), (\cC\! : \!A))$-module-morphisms} $\bSi\to \tbSi$
 is defined as a
$(\cD,\cC)$-bicomodule map $\cD\ot_B\Sigma\to \cD\ot_B\tSi$. The
 category of  $((\cD\! : \!B), (\cC\! : \!A))$-module-morphisms
 (with arrows given by module-morphism maps
and composed as bicolinear maps) is denoted by ${}_{(\cD :
B)}\cM_{(\cC : A)}$.

\begin{proposition}\label{pro.first-isom}
The category
of $((\cD\! : \!B), (\cC\! : \!A))$-module-morphisms ${}_{(\cD :
B)}\cM_{(\cC : A)}$ is isomorphic to the category
${}_{(\cD : B)}{\mathcal{R}^{\mathsf{f}}}_{(\cC : A)}$.
\end{proposition}
\begin{proof}
On objects, mutually inverse functors are given by the bijective
correspondence of Lemma~\ref{lemma.mod-mor}.
The
actions on morphisms are given by the bijections of the hom-sets
stated in Proposition~\ref{red-unred}.
\end{proof}

The duality in Lemma~\ref{duality} and observations made in
Remark~\ref{left-1-cell} allow one to construct a category dual to
${}_{(\cD :
B)}\cM_{(\cC : A)}$.
The resulting category
generalises that
of representations of a coring in a coring (cf.\ \cite[24.3]{BrzWis:cor}).
By Remark~\ref{left-1-cell}, any module-morphism
${}_{(\cD : B)}\bSi_{(\cC : A)}$ can be viewed as an object in
${}_{(\cD : B)}{\mathcal{L}^{\mathsf{f}}}_{(\cC : A)}$. The hom-set bijections
\begin{equation}\label{eq.2-iso}
\hom{A,A}{\Sigma^*\ot_B\cD\ot_B\tSi}{A} \,\, \simeq \,\,
\hom{A,B}{\Sigma^*\ot_B\cD}{{\tSi}^*},
\end{equation}
allow for identification of morphisms in
${}_{(\cD : B)}{\mathcal{L}^{\mathsf{f}}}_{(\cC : A)}$ with
$A$-bimodule maps
\begin{equation}
f: \Sigma^*\ot_B \cD \ot_B \tSi \to A,
\label{eq.f}
\end{equation}
such that for all $s^*\in \Sigma^*$, $\tilde{s}\in \tSi$ and
$d\in\cD$,
\begin{equation}
\label{cond.mor} \sum_j f(s^*\ot_B d\sw 1\ot_B
\tilde{e}_j)\tsi(\tilde{e}_j^*\ot_B d\ot_B \tilde{s})  = \sum_i
\sigma(s^*\ot_Bd\sw 1\ot_B e_i)f(e^*_i\ot_B d\sw 2\ot_B
\tilde{s}).
\end{equation}
where
$\{e_i,e^*_i\}$ is a dual basis of $\Sigma$ and
$\{\tilde{e}_j,\tilde{e}^*_j\}$ is a dual basis of $\tSi$.
The
 com\-po\-si\-tion of $f$ in equation \eqref{eq.f} with
 $g: {\widetilde{\Sigma}}^*\tensor{B}\cD\tensor{B}\widehat{\Sigma}\to A$,
transferred from that in
 ${}_{(\cD : B)}{\mathcal{L}^{\mathsf{f}}}_{(\cC : A)},$
comes out as
\begin{equation}
\label{comp} g\diamond f: \Sigma^*\ot_B\cD\ot_B \hSi\to A, \quad
s^*\ot_Bd\ot_B \hat{s}\mapsto \sum_j f(s^*\ot_B d\sw 1\ot_B
\te_j)g(\te^*_j\ot_Bd\sw 2\ot_B \hat{s}).
\end{equation}
In view of Lemma~\ref{duality} and Proposition~\ref{pro.first-isom}, the
category $\MRep (\cD\! : \!B |\cC\! : \!A)$ with objects module-morphisms and
 arrows given as maps $f$ satisfying condition \eqref{cond.mor} and composed
 according to \eqref{comp}, is dual (in fact, anti-isomorphic) to
 ${}_{(\cD :
B)}\cM_{(\cC : A)}$.
Since any morphism can be viewed as a module-morphism
as in Example~\ref{ex.mod-mor}, the category $\Rep (\cD\! : \!B |\cC\! : \!A)$ of representations of
$(\cD\! : \!B)$ in $(\cC\! : \!A)$ (cf.\ \cite[24.3]{BrzWis:cor})
is a full subcategory of $\MRep (\cD\! : \!B |\cC\! : \!A)$.

\subsection{Push-out and pull-back functors}\label{out-back}

Since any module-morphism ${}_{(\cD : B)}\bSi_{(\cC : A)}$ induces a $1$-cell
in the
bicategory $\fREM$ by Lemma~\ref{lemma.mod-mor},
Proposition~\ref{push-Out} yields the existence of the functor
$\bSi_{\circ}: \rcomod{\cD} \rightarrow \rcomod{\cC}$,
$(M,\varrho^M) \mapsto \left(M\ot_B\Sigma, \varrho^{M\ot_B\Sigma}
\right)$. In terms of $\sigma$, the coaction is derived from equation
\eqref{fks-sigma}, and reads
\begin{equation}\label{eq.right}
\varrho^{\bSi_\circ (M)}: M\ot_B\Sigma\to M\ot_B\Sigma\ot_A\cC, \;
m\ot_B s \mapsto \sum_i m\sw 0\ot_B e_i\ot_A\sigma(e^*_i\ot_Bm\sw
1\ot_B s).
\end{equation}
The functor
$\bSi_\circ $ is called a {\em (right) push-out functor}. Taking
into account that $\sigma$ is a morphism of corings and the observations made in
Remark~\ref{coendomorph}, it is reasonable to compare the cotensor
functor induced by the quasi-finite bicomodule $\cD\ot_B\Sigma$
and the right push-out functor $\bSi_\circ$. This leads to
 the following

\begin{proposition}\label{coten=pushout}
If ${}_{(\cD : B)}\bSi_{(\cC : A)}$ is a module-morphism, then
$\cD\ot_B\Sigma$ induces a functor
$-\cotensor{\cD}(\cD\ot_B\Sigma): \rcomod{\cD} \to \rcomod{\cC}$
which is naturally isomorphic to  $
\bSi_\circ$.
\end{proposition}
\begin{proof}
Start with an arbitrary right $\cD$-comodule $M$ and
consider the following commutative diagram
$$ \xymatrix@R=30pt@C=50pt{ 0 \ar@{->}[r] & M \Box_\cD (\cD\ot_B\Sigma)
\ar@{->}^-{\omega_{M,\cD\ot_B\Sigma}^k}[r] & M
\ot_B(\cD\ot_B\Sigma) \ar@{->}^-{\omega_{M,\cD\ot_B\Sigma}}[r]
\ar@<0,5ex>[dl]^-{M\ot_B\varepsilon_{\cD}\ot_B\Sigma} & M \ot_B \cD\ot_B(\cD\ot_B\Sigma) \\
& M\ot_{B}\Sigma \ar@{-->}^-{f}[u]
\ar@<0,5ex>[ur]^-{\varrho^M\ot_{B}\Sigma} & & }
$$
where $f$ is uniquely determined by the universal property of the
kernel $\omega_{M,\cD\ot_B\Sigma}^k$ in the category of right
$A$-modules. Clearly, $(M\ot_B\eD\ot_B\Sigma) \circ
\omega_{M,\cD\ot_B\Sigma}^k$ is the inverse map of $f$, and so $f$
is an isomorphism of right $A$-modules. Consider the composition
map
$$
\xymatrix@R=30pt@C=50pt{
M\Box_\cD(\cD\ot_B\Sigma)\ar@{->}_-{f^{-1}}[d]
\ar@{-->}^-{\varrho'}[r] &  \left( M\Box_\cD(\cD\ot_B\Sigma) \right)\ot_A\Sigma[\cD] \\
M\ot_B\Sigma \ar@{->}^-{\eta_M\ot_B\Sigma}[r] &
(M\ot_B\Sigma)\ot_A\Sigma[\cD] \ar@{->}_-{f\ot_A\Sigma[\cD]}[u]  }
$$
where $\eta_{-}$ is the natural transformation in equation
(\ref{eq.unit}) defined at $M$ by $$\eta_M(m\ot_B
s)=\sum_{i}m_{(0)}\ot_B e_i\ot_A e_i^*\ot_B m_{(1)},$$ for
every $m\in M$, $s \in \Sigma$.
Since the coaction $\varrho^{\cD\ot_B\Sigma}$
of Lemma \ref{comodule} satisfies $\eta_{M\ot_B\cD}\ot_B\Sigma=M
\ot_B\varrho^{\cD\ot_B\Sigma}$, we obtain the following commutative diagram
$$
\xymatrix@R=50pt@C=40pt{ & M\Box_\cD(\cD\ot_B\Sigma)
\ar@{->}^-{\varrho'}[r] \ar@{_{(}->}[dl] &
\left(M\Box_\cD(\cD\ot_B\Sigma)\right)\ot_A\Sigma[\cD] \ar@{_{(}->}[dl] \\
M\ot_B(\cD\ot_B\Sigma)
\ar@{->}^-{M\ot_B\varrho^{\cD\ot_B\Sigma}}[r] &
(M\ot_B(\cD\ot_B\Sigma))\ot_A\Sigma[\cD] & \\ & M\ot_B\Sigma
\ar@{_{(}->}^-{\varrho^M\ot_B\Sigma}[ul]
\ar@{->}^-{\eta_M\ot_B\Sigma}[r] \ar@{->}'[u]^>>>>{f} [uu] &
(M\ot_B\Sigma)\ot_A\Sigma[\cD]
\ar@{_{(}->}_-{\varrho^M\ot_B\Sigma\ot_A\Sigma[\cD]}[ul]
\ar@{->}_-{f\ot_A\Sigma[\cD]}[uu] }
$$
The right $\Sigma[\cD]$-coaction $M\ot_B\varrho^{\cD\ot_B\Sigma}$
induces the structure of a right $\Sigma[\cD]$-comodule on both
$M\Box_\cD(\cD\ot_B\Sigma)$ and $M\ot_B\Sigma$. The coactions are
given, respectively, by $\varrho'$ and $\eta_M\ot_B\Sigma$.
Furthermore, the $A$-module isomorphism $M\Box_\cD(\cD\ot_B\Sigma)
\simeq M\ot_B\Sigma$ becomes an isomorphism of right
$\Sigma[\cD]$-comodules. Applying the induction functor associated
to the underlying morphism of corings $\sigma$ to this
$\Sigma[\cD]$-comodule isomorphism, we obtain an isomorphism of
right $\cC$-comodules which gives the desired natural
isomorphism.
\end{proof}

By Remark~\ref{left-1-cell} and the left-handed
versions of both Proposition~\ref{push-Out} and
Proposition~\ref{coten=pushout}, a module-morphism
 ${}_{(\cD : B)}\bSi_{(\cC : A)}$ induces
a {\em (left) push-out functor}
$$
{}_\circ\bSi:{}^\cD\cM\to {}^\cC\cM, \qquad N\mapsto
(\Sigma^*\ot_B\cD)\Box_\cD N \simeq \Sigma^*\ot_BN.
$$
Explicitly, using equation
\eqref{fkt}, the left $\cC$-coaction reads, for all $s^*\in \Sigma^*$
and $n\in N$,
\begin{equation}\label{eq.left.coa}
{}^{{}_\circ\bSi(N)}\varrho (s^*\ot_B n) = \sum_i \sigma(s^*\ot_B n\sw {-1}\ot_B e_i)
\ot_A e^*_i\ot_B n\sw 0.
\end{equation}
In particular, in view of Lemma~\ref{comodule},
${}_\circ\bSi(\cD)$ is a $(\cC,\cD)$-bicomodule. We say that a module-morphism
${}_{(\cD : B)}\bSi_{(\cC : A)}$ is {\em right pure} if for all right $\cC$-comodules $L$
the coaction equalising map
$$\omega_{L,{}_\circ\bSi(\cD)} =\varrho^L \ut_A {}_\circ\bSi(\cD) - L
\ut_A {}^{{}_\circ\bSi(\cD)}\!\varrho : L\ot_A {}_\circ\bSi(\cD) \to
L\ot_A\cC\ot_A {}_\circ\bSi(\cD)
$$
is a $\cD\ot_B\cD$-pure morphism of right $B$-modules. Obviously,
if $\cD$ is a flat left $B$-module, then every right $B$-module
map is $\cD\ot_B\cD$-pure, hence any module-morphism ${}_{(\cD :
B)}\bSi_{(\cC : A)}$ is right pure in this case.

If ${}_{(\cD : B)}\bSi_{(\cC : A)}$ is  right pure then for all
right $\cC$-comodules $N$, $N\Box_\cC {}_\circ\bSi(\cD)$ is a
right $\cD$-comodule (cf.\ \cite[22.3]{BrzWis:cor}). Thus any
right pure module-morphism ${}_{(\cD : B)}\bSi_{(\cC : A)}$ gives
rise to a functor
 $$
 \bSi^\circ  : \cM^\cC\to \cM^\cD, \qquad N\mapsto N\Box_\cC {}_\circ\bSi(\cD).
 $$
The $\cD$-coaction on $\bSi^\circ (N)$, $ \varrho^{\bSi^\circ
(N)}: N\Box_\cC  ({}\Sigma^*\ot_B\cD)\to (N\Box_\cC
({}\Sigma^*\ot_B\cD))\ot_B\cD$ explicitly reads:
 $$
\varrho^{\bSi^\circ (N)}(\sum_\alpha n^\alpha\ot_A s^*_\alpha\ot_B d^\alpha) \sum_\alpha n^\alpha\ot_A s^*_\alpha\ot_B d^\alpha\sw 1\ot_B d^\alpha\sw 2.
$$
The functor $\bSi^\circ $ is called a {\em (right) pull-back functor}
associated to a right pure module-morphism ${}_{(\cD : B)}\bSi_{(\cC : A)}$.

In a similar way one defines a left pure morphism and a (left) pull-back functor. There is
an obvious left-right symmetry, hence we restrict ourselves to right pure morphisms
and (right) pull-back functors.

\begin{theorem}
\label{thm.adjoint}
For any right pure module-morphism ${}_{(\cD : B)}\bSi_{(\cC : A)}$,
the pull-back functor $\bSi^\circ $ is the right adjoint to the push-out functor
$\bSi_\circ $.
\end{theorem}
\begin{proof}
First we construct the unit of the adjunction. For any $M\in\cM^\cD$,
consider a $k$-linear map
$$
\eta_M :M\to \bSi^\circ (\bSi_\circ (M)) (M\ot_B\Sigma)\Box_\cC(\Sigma^*\ot_B\cD),\quad m\mapsto
\sum_i m\sw 0\ut_B e_i\ot_Ae^*_i\ut_B m\sw 1.
$$
Clearly, the map $\eta_M$ is
well-defined and it immediately follows from the definitions of
the coactions on $\bSi_\circ (M)$ and ${}_\circ\bSi (\cD)$ that the image
of $\eta_M$ is in the required cotensor product.  The way in which
the definition of the map $\eta_ M$ depends
upon the coaction $\varrho^M$ ensures that $\eta_M$ is a right
$B$-module map (since $\varrho^ M$ is such a map). One easily checks that
$\eta_ M$ is also a morphism in $\cM^\cD$.
 Next, take any $f:  M\to  M'$ in $\cM^\cD$ and compute for any
$m\in  M$, $$\begin{array}{rcl}
\eta_{ M'}( f(m))& = & \sum_i f(m)\sw 0\ut_B e_i\ot_Ae^*_i\ut_B f(m)\sw 1 \\[+1mm]
   &  & \sum f(m\sw 0)\ut_B e_i\ot_Ae^*_i\ut_B m\sw 1\; =\;
\bSi_\circ (\bSi^\circ (f))(\eta_ M(m)).
\end{array}$$
The second equality follows, since $f$ is a morphism of right
$\cD$-comodules. Thus we have constructed a natural map $\eta:
\id_{\cM^\cD}\to \bSi_\circ  \bSi^\circ $ that will be shown to be
the unit of the adjunction.

Now, for any $ N\in \cM^{\cC}$, consider a right $A$-module map
$$
\begin{array}{c}
\psi_{ N}: \bSi_\circ (\bSi^\circ ( N))=(N\Box_\cC(\Sigma^*\ot_B\cD))\ot_B\Sigma \to  N,
\\ [+2mm]
{\sum}_{\alpha}n^{\alpha}\otimes_{A}s^*_\alpha\ot_Bd^{\alpha}\otimes_{B}s\mapsto
{\sum}_\alpha n^{\alpha}s^*_\alpha(\eD(d^{\alpha})s).
\end{array}$$
First we need to check whether the map $\psi_{ N}$ is a morphism in the
category of right $\cC$-comodules. Take any
$x= {\sum}_{\alpha}n^{\alpha}\otimes_{A}s^*_\alpha\ot_Bd^{\alpha}\otimes_{B}s
\in (N\Box_\cC(\Sigma^*\ot_B\cD))\ot_B\Sigma$. Then,
\begin{eqnarray*}
\sum\psi_N(x\sw 0)\ot_B x\sw 1 &=&  \sum_{i,\alpha} \psi_{ N}(n^{\alpha}
\otimes_{A}s^*_\alpha\ot_Bd^{\alpha}\sw 1\otimes_{B}e_i)\ot_A\sigma(e^*_i
\ot_B d^\alpha\sw 2\ot_B s)
\\ &=&     \sum_{i,\alpha} n^{\alpha}s^*_\alpha(\eD(d^{\alpha}\sw 1)e_i)
\ot_A\sigma(e^*_i\ot_B d^\alpha\sw 2\ot_B s)\\
   &=&  \sum_{i,\alpha} n^{\alpha}\ot_A s^*_\alpha(e_i)
   \sigma(e^*_i\ot_B d^\alpha\ot_B s) =  \sum_{\alpha} n^{\alpha}\ot_A
   \sigma(s^*_\alpha\ot_B d^\alpha\ot_B s).
\end{eqnarray*}
The final equality is a consequence of the left $A$-linearity of $\sigma$ and the dual
basis property.
Since ${\sum}_{\alpha}n^{\alpha}\otimes_{A}s^*_\alpha\ot_Bd^{\alpha}
\in N\Box_\cC(\Sigma^*\ot_B\cD)$,
$$
{\sum}_{\alpha}n^{\alpha}\sw 0\ot_A n^\alpha\sw 1\otimes_{A}s^*_\alpha\ot_Bd^{\alpha}
 = {\sum}_{\alpha,i}n^{\alpha}\otimes_{A}
 \sigma(s^*_\alpha\ot_Bd^{\alpha}\sw 1\ot_B e_i)\ot_A e^*_i\ot_B d^\alpha\sw 2,
$$
hence
$$
{\sum}_{\alpha}n^{\alpha}\sw 0\ot_A n^\alpha\sw 1\otimes_{A}s^*_\alpha\eD(d^{\alpha})
 = {\sum}_{\alpha,i}n^{\alpha}\otimes_{A}\sigma(s^*_\alpha\ot_Bd^{\alpha}\ot_B e_i)
 \ot_A e^*_i.
$$
Using this equality, the $A$-linearity of $\sigma$ and the properties
of a dual basis, we can compute
\begin{eqnarray*}
   \sum\psi_N(x\sw 0)\ot_B x\sw 1 &=&   \sum_{\alpha} n^{\alpha}\ot_A
   \sigma(s^*_\alpha\ot_B d^\alpha\ot_B s)
   = \sum_{\alpha,i} n^{\alpha}\ot_A \sigma(s^*_\alpha\ot_B d^\alpha\ot_B e_i e^*_i(s))\\
   &=& \sum_{\alpha,i} n^{\alpha}\ot_A \sigma(s^*_\alpha\ot_B d^\alpha\ot_B e_i) e^*_i(s)\\
   &=& {\sum}_{\alpha}n^{\alpha}\sw 0\ot_A n^\alpha\sw 1s^*_\alpha(\eD(d^{\alpha})s)
   = \sum\psi_N(x)\sw 0\ot_A \psi_N(x)\sw 1 .
\end{eqnarray*}
Therefore, $\psi_N$ is a right $\cC$-comodule map as required. Thus for any right
$\cC$-comodule $ N$  we have
constructed a morphism $\psi_{ N}$ in $\cM^{\cC}$.
Noting that any right $\cC$-comodule map is necessarily a right $A$-module map,
one easily checks that the map $\psi_{ N}$ is natural in $N$, i.e.\
the collection of all the $\psi_{ N}$ defines a morphism of functors
$\psi : \bSi_\circ \bSi^\circ \to \id_{\cM^{\cC}}$.

The verification that $\eta$ and $\psi$ are the unit and
counit  respectively, i.e.\ that for all $M\in \cM^\cD$ and $N\in \cM^\cC$,
$
\psi_{\bSi_\circ (M)}\circ \bSi_\circ (\eta_M) = \bSi_\circ (M)$ and
$\bSi^\circ (\psi_N)\circ \eta_{\bSi^\circ (N)} =\bSi^\circ (N)$,
is a straightforward application of the properties of a dual basis, and is left
to the reader.
\end{proof}

\begin{remark}
One easily checks that the unit and counit defined in the
proof of Theorem~\ref{thm.adjoint} are induced, respectively, by
the unit and counit of the adjunction (\ref{eq.adj1}) (i.e.\ by
the natural transformations (\ref{eq.unit}) and
(\ref{eq.counit})). If a module-morphism
${}_{(\cD : B)}\bSi_{(\cC : A)}$ is right pure, then
Theorem~\ref{thm.adjoint} asserts that this last adjunction is
extended to the category of right $\cC$-comodules. That is,
Theorem~\ref{thm.adjoint} and its proof can be seen as an example and
also an application of the statements and the proof of
\cite[Proposition 4.2(1)]{Gom:sep}.
\end{remark}

In case ${}_{(\cD : B)}\bSi_{(\cC : A)}$ corresponds to a morphism
of corings, i.e.\ there is an algebra map $B\to A$, $\Sigma = A$
and $\sigma: A\ot_B \cD\ot_B A\to \cC$ is an $A$-coring map (cf.\
Example~\ref{ex.mod-mor}), the functor $\bSi_\circ $ is the
induction functor and $\bSi^\circ $ is the ad-induction functor
introduced in \cite{Gom:sep}. In this case
Theorem~\ref{thm.adjoint} reduces to \cite[Proposition
5.4]{Gom:sep}.

Theorem~\ref{thm.adjoint} tantamounts  to the existence of
isomorphisms, for all $M\in \cM^\cD$ and $N\in \cM^\cC$,
\begin{equation}\label{omega1}
\Omega_{M,N} : \Rhom \cC {M\ot_B \Sigma}N\to \Rhom \cD M {N\Box_\cC( \Sigma^*\ot_B \cD)},
\end{equation}
for a  right pure module-morphism ${}_{(\cD : B)}\bSi_{(\cC : A)}$.
Explicitly, these isomorphisms read, for all $\phi\in \Rhom \cC {M\ot_B \Sigma}N$
and  $m\in M$,
\begin{equation}\label{omega2}
\Omega_{M,N}(\phi) (m) = \sum_i \phi(m\sw 0\ot_B e_i)\ot_A e^*_i\ot_B m\sw 1.
\end{equation}
To write out the inverse of $\Omega_{M,N}$ explicitly,
take any $\tilde{\phi}\in \Rhom \cD M {N\Box_\cC( \Sigma^*\ot_B \cD)}$
and $m\in M$, and write $\tilde{\phi}(m) \sum \tilde{\phi}(m)\su 1\ot_A \tilde{\phi}(m)\su 2\ot_B \tilde{\phi}(m)\su 3$.
Then, for all $s\in \Sigma$,
\begin{equation}\label{omega3}
\Omega_{M,N}^{-1}(\tilde{\phi}) (m\ot_B s) = \sum
\tilde{\phi}(m)\su 1\tilde{\phi}(m)\su 2(\eD(\tilde{\phi}(m)\su 3)s).
\end{equation}
A module-morphism ${}_{(\cD : B)}\bSi_{(\cC : A)}$ is assumed to be
right pure to assure that for all right $\cC$-comodules $N$,
$\bSi^\circ (N)$ is a right $\cD$-comodule. On the other hand,
for any module-morphism there exist right $\cC$-comodules $N$,
such that $\bSi^\circ (N)$ is a right $\cD$-comodule.
For any such comodule the isomorphisms in equations (\ref{omega1})-(\ref{omega3})
are well-defined. This observation leads to the following

\begin{corollary}
\label{cor.adjoint}
For any  module-morphism ${}_{(\cD : B)}\bSi_{(\cC : A)}$ and any right $\cD$-comodule $M$,
$$
\Rhom \cD M{\Sigma^*\ot_B\cD} \simeq \Rhom \cC{M\ot_B \Sigma}
\cC \simeq \rhom A {M\ot_B \Sigma} A =(M\ot_B \Sigma)^* ,
$$
as $k$-modules.
\end{corollary}
\begin{proof}
Note that $\bSi^\circ (\cC) = \cC\Box_\cC (\Sigma^*\ot_B\cD)\simeq \Sigma^*\ot_B\cD$,
so that it is a right $\cD$-comodule. Thus $\Omega_{M,\cC}:
\Rhom \cC{M\ot_B \Sigma} \cC \to \Rhom \cD M{\Sigma^*\ot_B\cD}$
are well-defined isomorphisms. In view of equations (\ref{omega2})-(\ref{omega3}),
and the identification $ \cC\Box_\cC (\Sigma^*\ot_B\cD)\simeq \Sigma^*\ot_B\cD$,
these isomorphisms and their inverses come out explicitly as, for all $m\in M$
and $\phi \in \Rhom \cD M{\Sigma^*\ot_B\cD}$,
$$
\Omega_{M,\cC} (\phi)(m) =  \sum_i \eC(\phi(m\sw 0\ot_B e_i)) e^*_i\ot_B m\sw 1,
$$
and for all $\tilde{\phi}\in \Rhom \cC{M\ot_B \Sigma} \cC$, $m\in M$ and $s\in \Sigma$,
$$
\Omega_{M,\cC}^{-1}(\tilde{\phi}) (m\ot_B s) = \sigma(\tilde{\phi}(m)\ot_B s).
$$
The second isomorphism follows from $\Rhom\cC N\cC\simeq \rhom A N
A$, which holds for any right $\cC$-comodule $N$ (and is given by
the counit with the inverse provided by  the coaction).
\end{proof}

\begin{corollary}\label{cor.ring}
For any  module-morphism ${}_{(\cD : B)}\bSi_{(\cC : A)}$,
the endomorphism ring of right $\cD$-comodule $\Sigma^*\ot_B\cD$ is
isomorphic to the right dual ring of the $A$-coring $\Sigma[\cD]$, i.e.\
$
\Rend \cD {\Sigma^*\ot_B\cD} \simeq (\Sigma^*\ot_B\cD\ot_B\Sigma)^*,
$
as rings.
\end{corollary}
\begin{proof}
Setting $M= \Sigma^*\ot_B\cD$ in Corollary~\ref{cor.adjoint},
one obtains an isomorphism of $k$-modules
$$
\Theta: \Rend \cD {\Sigma^*\ot_B\cD} \to (\Sigma^*\ot_B\cD\ot_B\Sigma)^*,
\;\; {\phi} \mapsto [ s^*\ot_B d\ot_B s \mapsto \eC(\sigma(\phi(s^*\ot_B d)\ot_B s))].
$$
We need to check if this isomorphism is an algebra map. For any
$\phi\in \Rend \cD {\Sigma^*\ot_B\cD}$, $d\in \cD$ and $s^*\in \Sigma^*$,
write $\phi(s^*\ot_B d) = \sum \phi(s^*\ot_B d)\su 1\ot_B \phi(s^*\ot_B d)
\su 2 \in \Sigma^*\ot_B \cD$. Since $\sigma$ is a morphism of $A$-corings,
$$
\eC(\sigma(\phi(s^*\ot_B d)\ot_B s)) = \eSD(\phi(s^*\ot_B d)\ot_B s) \sum \phi(s^*\ot_B d)\su 1(\eD( \phi(s^*\ot_B d)\su 2s)).
$$
In particular $\Theta$ maps the identity morphism in $\Sigma^*\ot_B\cD$ to
the counit $\eSD$. Furthermore, since $\phi$ is a right $\cD$-comodule map
$$
\phi(s^*\ot_Bd\sw 1)\ot_B d\sw 2 = \sum \phi(s^*\ot_B d)\su 1\ot_B
\phi(s^*\ot_B d)\su 2 \sw 1\ot_B \phi(s^*\ot_B d)\su 2 \sw 2,
$$
so that
$$
 \sum \phi(s^*\ot_B d\sw 1)\su 1\eD( \phi(s^*\ot_B d\sw 1)\su 2) \ot_B d\sw 2 = \phi(s^*\ot_B d).
 $$
 Taking this into account we can compute, for all $\phi,\phi'\in \Rend \cD {\Sigma^*\ot_B\cD}$,
 $d\in \cD$, $s\in \Sigma$ and $s^*\in \Sigma^*$,
 \begin{eqnarray*}
( \Theta(\phi')*\Theta(\phi))(s^*\ot_B d\ot_B s) &=& \sum_i \Theta(\phi')(\Theta(\phi)(s^*\ot_B d\sw 1\ot_B e_i)e_i^*\ot_B d\sw 2\ot_B s)\\
&\hspace{-1.1in}=&\hspace{-.55in} \sum_i \Theta(\phi')(\phi(s^*\ot_B d\sw 1)\su 1(\eD(\phi(s^*\ot_B d\sw 1)\su 2) e_i)e_i^*\ot_B d\sw 2\ot_B s)\\
&\hspace{-1.1in}=&\hspace{-.55in} \sum_i \Theta(\phi')(\phi(s^*\ot_B d)\su 1( e_i)e_i^*\ot_B \phi(s^*\ot_B d)\su 2\ot_B s)\\
&\hspace{-1.1in}=&\hspace{-.55in} \sum \Theta(\phi')(\phi(s^*\ot_B d)\su 1\ot_B \phi(s^*\ot_B d)\su 2\ot_B s)\\
&\hspace{-1.1in}=&\hspace{-.55in} \sum \Theta(\phi')(\phi(s^*\ot_B d)\ot_B s) = \eC(\sigma(\phi'(\phi(s^*\ot_B d))\ot_B s))\\
&\hspace{-1.1in}=&\hspace{-.55in} \Theta(\phi'\circ\phi)(s^*\ot_B d\ot_B s).
\end{eqnarray*}
Thus we conclude that $\Theta$ is an algebra map as required.
\end{proof}

In case $\cD=B$, Corollary~\ref{cor.ring} implies that the ring of
endomorphisms of the right $B$-module $\Sigma^*$ is isomorphic to
the right dual ring of the corresponding comatrix coring.
This is a right-handed version of one of the assertions in
\cite[Proposition~2.1]{KaoGom:com}.

\subsection{Natural isomorphisms between push-out
functors}\label{subsec.nat-iso} The horizontal composition in
bicategory $\fREM$ induces the {\em cotensor product} of
 module-morphisms. Start with morphisms $\phi
:\bSi\to\tbSi$ in ${}_{(\cD : B)}\cM_{(\cC : A)}$ and $\psi:
\bXi\to \tbXi$ in ${}_{(\cE : C)}\cM_{(\cD : B)}$, i.e.\ $\phi :
\cD\ot_B\Sigma\to \cD\ot_B \tSi$ is a $(\cD,\cC)$-bicomodule map
and $\psi : \cE\ot_C\Xi\to \cE\ot_C \tXi$ is a
$(\cE,\cD)$-bicomodule map. Therefore, we can consider the
$(\cE,\cC)$-bicomodule map
$$
\psi\Box_\cD\phi: (\cE\ot_C\Xi)\Box_\cD(\cD\ot_B\Sigma)\to (\cE\ot_C\tXi)\Box_\cD(\cD\ot_B\tSi).
$$
Using the isomorphisms $\cE\ot_C\Xi\ot_B\Sigma\simeq
(\cE\ot_C\Xi)\Box_\cD(\cD\ot_B\Sigma)$ and
$(\cE\ot_C\tXi)\Box_\cD(\cD\ot_B\tSi)\simeq \cE\ot_C\tXi\ot_B\tSi$
we thus arrive at an $(\cE,\cC)$-bicomodule map, hence the
map of module-morphisms,
$$
\psi\Box_\cD\phi: \cE\ot_C\Xi\ot_B\Sigma\to \cE\ot_C\tXi\ot_B\tSi.
$$

The morphisms in the category of module-morphisms detect
relationships between push-out functors. More precisely, one can
formulate the following

\begin{proposition}\label{prop.mor.nat}
There is a bijective correspondence between morphisms in the
category of module-morphisms  and natural transformations between
corresponding push-out functors. This correspondence is compatible
with horizontal compositions.
\end{proposition}
\begin{proof}
Let ${}_{(\cD : B)}\bSi_{(\cC : A)}$, ${}_{(\cD : B)}\tbSi_{(\cC :
A)}$ be module-morphisms and suppose that $f: \bSi_\circ
\to\tbSi_\circ $ is a natural map. This means that for all $M\in
\cM^\cD$, there is a right $\cC$-comodule map $f_M:M\ot_B\Sigma\to
M\ot_B\tSi$. Consider morphisms in $\cM^\cD$, $\varrho^M:M\to
M\ot_B \cD$ and, for all $m\in M$, $\ell_m :\cD\to M\ot_B \cD$,
$d\mapsto m\ot_B d$. The naturality of $f$ implies that, for all
$m\in M$, $d\in \cD$ and $s\in\Sigma$, $$
f_{M\ot_B\cD}(m\ot_Bd\ot_Bs) = m\ot_B f_\cD(d\ot_B s), \quad \sum
f_{M\ot_B\cD}\circ (\varrho^M \ot_B\Sigma) (\varrho^M\ot_B\tSi)\circ f_M.$$ Put together, this means that
\begin{equation}
(\varrho^M\ot_B\tSi)\circ f_M = (M\ot_B f_\cD)\circ
(\varrho^M\ot_B \Sigma).
\label{eq.nat}
\end{equation} If $M=\cD$ equation \eqref{eq.nat}
implies that $f_\cD$ is a left $\cD$-comodule map, hence it is a
$(\cD,\cC)$-bicomodule map. Recall that
$M\Box_\cD(\cD\ot_B\Sigma)\simeq M\ot_B\Sigma$ by $M\ot_B \eD\ot_B
\Sigma$ and $\varrho^M\ot_B \Sigma$. Thus, applying $M\ot_B
\eD\ot_B \tSi$ to equation \eqref{eq.nat} we obtain
$$
 f_M= (M\ot_B \eD\ot_B \tSi)\circ (M\ot_B f_\cD)\circ (\varrho^M\ot_B \Sigma) \simeq M\Box_\cD f_\cD.
$$
Hence the required bijective correspondence is provided by
$$
{\rm Nat}(\bSi_\circ ,\tbSi_\circ ) \ni f \mapsto f_\cD \in \LRhom
\cD \cC{\cD\ot_B\tSi}{\cD\ot_B\Sigma},
$$
with the inverse, for all $M\in \cM^\cD$, $$ \LRhom
\cD \cC{\cD\ot_B\tSi}{\cD\ot_B\Sigma}\ni\phi\mapsto (M\ot_B \eD\ot_B
\tSi)\circ (M\ot_B \phi)\circ (\varrho^M\ot_B \Sigma) \simeq
M\Box_\cD\phi.
$$
Note that this is natural in $M$ by the functoriality of the
cotensor product. Clearly, this bijective correspondence is
compatible with compositions.
\end{proof}

\begin{corollary}\label{cor.mor.nat}
Let ${}_{(\cD : B)}\bSi_{(\cC : A)}$, ${}_{(\cD : B)}\tbSi_{(\cC :
A)}$ be  $((\cD\! :\! B),(\cC\! :\! A))$-module-morphisms. Then the following
statements are equivalent:
\begin{blist}
\item push-out functors $\bSi_\circ $ and $\tbSi_\circ $ are
naturally isomorphic to each other;
\item  ${}_{(\cD :
B)}\bSi_{(\cC : A)}\simeq {}_{(\cD : B)}\tbSi_{(\cC : A)}$  in
$\MRep (\cD\! : \!B |\cC\! : \!A)$;
\item  ${}_{(\cD :
B)}\bSi_{(\cC : A)}\simeq {}_{(\cD : B)}\tbSi_{(\cC : A)}$  in $
{}_{(\cD : B)}\cM_{(\cC : A)}$.
\end{blist}
\end{corollary}
\begin{proof}
This follows immediately from Proposition~\ref{prop.mor.nat} and
the duality between $\MRep (\cD\! : \!B |\cC\! : \!A)$ and $
{}_{(\cD : B)}\cM_{(\cC : A)}$.
\end{proof}

\section{Equivalences induced by push-out and pull-back functors}\label{sec.equiv}
In this section we study when a push-out functor is an equivalence. In terms of
 non-commutative algebraic geometry, this is a problem of determining which
 changes of covers are admissible (a change of
cover of a non-commutative space should not change the space, i.e.\ the associated
category of sheaves). We then proceed to study the generalised descent
associated to a morphism of corings.

\subsection{Criteria for an equivalence}
\begin{theorem}\label{thm.full}
Let ${}_{}\bSi = (\Sigma,\sigma)$ be a $((\cD\! :\! B),(\cC\! :\!
A))$-module-morphism  and assume that ${}_B\cD$ is flat.
\begin{zlist}
\item If $\sigma$ is an isomorphism of $A$-corings and
$\cD\ot_B\Sigma$ is a coflat left $\cD$-comodule, then  ${}_A\cC$
is flat and the pull-back functor $\bSi^\circ $ is full and
faithful. \item  If the pull-back functor $\bSi^\circ $ is full
and faithful then $\sigma$ is an isomorphism of $A$-corings.
\end{zlist}
\end{theorem}
\begin{proof}
First note that since ${}_B\cD$ is flat, $\bSi$ is a right pure
module-morphism so that the pull-back functor is well-defined.
Furthermore, $\cM^\cD$ is a Grothendieck category.

(1) If $\cD\ot_B\Sigma$ is a coflat left $\cD$-comodule, then the
cotensor functor $-\Box_\cD(\cD\ot_B\Sigma)$ is exact. Since
$M\Box_\cD (\cD\ot_B\Sigma)\simeq M\ot_B\Sigma$, every short exact
sequence
 $0\to M\to M'\to M''\to 0$ of right $\cD$-comodules yields an exact
 sequence  $0\to M\ot_B \Sigma\to M'\ot_B\Sigma \to M''\ot_B\Sigma\to 0$.

Consider an exact sequence $0\to V\to V'\to V''\to 0$ of right
$A$-modules. Since ${}_A\Sigma^*$ and ${}_B\cD$ are flat, the
above sequence yields an exact sequence $$ 0\to
V\ot_A\Sigma^*\ot_B\cD\to V'\ot_A\Sigma^*\ot_B\cD\to
V''\ot_A\Sigma^*\ot_B\cD\to 0
$$
of right $\cD$-comodules. The coflatness of $\cD\ot_B\Sigma$ then produces an exact sequence
$$
0\to V\ot_A\Sigma^*\ot_B\cD\ot_B\Sigma\to
V'\ot_A\Sigma^*\ot_B\cD\ot_B\Sigma\to
V''\ot_A\Sigma^*\ot_B\cD\ot_B\Sigma\to 0.
$$ Hence $\Sigma[\cD]$ is a flat left $A$-module, and since
$\sigma: \Sigma[\cD]\to \cC$ is an isomorphism of $A$-bimodules, also ${}_A\cC$ is flat.

For any right $\cC$-comodule $N$, consider the following
commutative diagram with exact rows
$$\xymatrix{
  0\ar[r] &(N\Box_\cC(\Sigma^*\ot_B\cD))\ot_B\Sigma \ar[r] \ar[d]^{\psi_N}  &
      N\ot_A \Sigma^*\ot_B\cD\ot_B\Sigma\ar[r] \ar[d]^{N\ot_A\sigma}  &
  N\ot_A \cC\ot_A \Sigma^*\ot_B\cD\ot_B\Sigma\ar[d]^{N\ot_A\cC\ot_A\sigma}   \\
  0\ar[r] & N \ar[r]^{\varrho^N} & N\ot_A\cC \ar[r]^{N\ot_A\DC - \varrho^N\ot_A\cC} & N\ot_A\cC\ot_A\cC \, .
} $$ The top row is the defining sequence of a cotensor product
tensored with ${}_B\Sigma$, hence it is exact (for
$\cD\ot_B\Sigma$ is coflat). The bottom row is exact by the
coassociativity of the coaction. Since $\sigma$ is an isomorphism
of right $\cC$-comodules, so are the second and third vertical
maps. This implies that the counit of the adjunction $\psi_N$ is
an isomorphism of right $\cC$-comodules. Thus the pull-back
functor $\bSi^\circ $ is full and faithful.

(2) If  $\bSi^\circ $ is a full and faithful functor, the counit
$\psi_N$ is an isomorphism of right $\cC$-comodules for any $N\in
\cM^\cC$. In particular, $\psi_\cC : (\cC\Box_\cC
(\Sigma^*\ot_B\cD))\ot_B\Sigma\to \cC$ is an isomorphism of right
$\cC$-comodules. Using the definitions of the left coaction
${}^{{}_\circ\bSi(\cC)}\varrho$ in equation (\ref{eq.left.coa})
and of $\psi_\cC$ from the proof of Theorem~\ref{thm.adjoint}, we
can compute, for all $s\in \Sigma$, $d\in\cD$ and
$s^*\in\Sigma^*$,
\begin{eqnarray*}
\psi_\cC\left(\left({}^{{}_\circ\bSi(\cC)}\varrho\ot_B\Sigma\right)\left(s^*\ot_Bd\ot_Bs\right)\right) &=&
\psi_\cC\left(\sum_i\sigma(s^*\ot_B d\sw 1\ot_B e_i)\ot_Ae^*_i\ot_Bd\sw 2\ot_B s\right)\\
&=& \sum_i\sigma(s^*\ot_B d\sw 1\ot_B e_i)e^*_i(\eD(d\sw 2)s)\\
&=& \sum_i\sigma(s^*\ot_B d\ot_B e_i)e^*_i(s) = \sigma(s^*\ot_B d\ot_Bs).
\end{eqnarray*}
Since ${}^{{}_\circ\bSi(\cC)}\varrho\ot_B\Sigma:
\Sigma^*\ot_B\cD\ot_B\Sigma \to (\cC\Box_\cC
(\Sigma^*\ot_B\cD))\ot_B\Sigma$ is an isomorphism, $\sigma$ is a
composition of isomorphisms, hence also an isomorphism (of $A$-corings). \end{proof}

\begin{theorem}\label{thm.equiv}
Let ${}_{}\bSi = (\Sigma,\sigma)$ be a $((\cD\! :\! B),(\cC\! :\!
A))$-module-morphism and assume that ${}_B\cD$ is flat. The
following statements are equivalent:
\begin{blist}
\item $\sigma$ is an isomorphism of $A$-corings and $\cD\ot_B\Sigma$ is a faithfully coflat left $\cD$-comodule;.
\item  ${}_A\cC$ is flat and  $\bSi^\circ $ is an equivalence of categories with the inverse $\bSi_\circ $.
\end{blist}
\end{theorem}
\begin{proof}
(a) $\Ra$ (b) By Theorem~\ref{thm.full}, $\cC$ is a flat left
$A$-module. We need to show that for all $M\in \cM^\cD$,
$N\in\cM^\cC$, a morphism $\phi\in \Rhom \cC{M\ot_B\Sigma}N$ is an
isomorphism if and only if $\Omega_{M,N}(\phi)\in \Rhom \cD
M{N\Box_\cC(\Sigma^*\ot_B\cD)}$ is an isomorphism. Here
$\Omega_{M,N}$ is the adjunction isomorphism given in
equation~(\ref{omega2}). Observe that there is an isomorphism
$$
\theta_N : N \to (N\Box_\cC(\Sigma^*\ot_B\cD))\Box_\cD(\cD\ot_B\Sigma),
$$
obtained as the following composition of isomorphisms:
 $$\xymatrix{
  N\ar[r]^{\varrho^N}&
  N\Box_\cC\cC\ar[r]^(.4){N\Box_\cC\sigma^{-1}}
  & N\Box_\cC(\Sigma^*\ot_B\cD\ot_B\Sigma) &\\
    ~~~\ar[rr]^(.3){N\Box_\cC(\Sigma^*\ot_B\DD\ot_B\Sigma)}& &
N\Box_\cC((\Sigma^*\ot_B\cD)\Box_\cD (\cD\ot_B\Sigma))\ar[r]&
(N\Box_\cC(\Sigma^*\ot_B\cD))\Box_\cD (\cD\ot_B\Sigma)  & .
  }
$$
The last isomorphism is a consequence of the fact that
$\cD\otimes_B\Sigma$ is a (faithfully) coflat left $\cD$-comodule.
In this way we are led to the following commutative diagram
  $$\xymatrix{  M \ot_B  \Sigma \ar[rr] ^(.4){\varrho^M\ot_B\Sigma}\ar[d]_{\phi}  & &
          M\square_\cD(\cD\otimes_B\Sigma)
  \ar[d]^{\Omega_{M,N}(\phi)\square_\cD(\cD\otimes_B\Sigma) } \\
      N \ar[rr] ^(.3){\theta_{N}}&
  &(N\Box_\cC(\Sigma^*\ot_B\cD))\Box_\cD (\cD\ot_B\Sigma) .
  } $$
  Since the rows are isomorphisms
  and $\cD\otimes_B\Sigma$ is a faithfully coflat left $\cD$-comodule, the map
  $\phi$ is an isomorphism if and only if $\Omega_{M,N}(\phi)$ is an isomorphism. Thus
  $\bSi^\circ $ is an equivalence as required.

  (b) $\Ra$ (a) Since ${}_A\cC$ and ${}_B\cD$ are flat, both $\cM^\cD$ and $\cM^\cC$
are Abelian categories, and kernels (and cokernels) are computed
in Abelian groups. The functor $\bSi_\circ -\Box_\cD(\cD\ot_B\Sigma)$ is an equivalence, hence it reflects
and preserves exact sequences. In view of the fact that a sequence
in $\cM^\cC$ is exact if and only if it is exact as a sequence of
Abelian groups, this implies that $\cD\ot_B\Sigma$ is a faithfully
coflat left $\cD$-comodule.
 \end{proof}

\begin{corollary}
\label{cor.equiv.com}
Let $\coring{D}$ be a $B$-coring, and $\Sigma$ a
$(B,A)$-bimodule such that $\Sigma_A$ is a finitely generated and
projective module.  If ${}_B\coring{D}$ is a flat module and
${}_B\Sigma$ is a faithfully flat module, then
\[
\xymatrix@C=30pt{ -\tensor{B}\Sigma: \rcomod{\coring{D}}
\ar@{->}[r] &  \rcomod{\Sigma[\coring{D}]} }
\]
is an equivalence of categories with the inverse
\[
\xymatrix@C=40pt{-\cotensor{\Sigma[\coring{D}]}(\Sigma^*\tensor{B}\coring{D}):
\rcomod{\Sigma[\coring{D}]} \ar@{->}[r] &
\rcomod{\coring{D}} .}
\]
In particular, if $B \rightarrow A$ is a left faithfully flat
ring extension (i.e.\ ${}_BA$ is a faithfully flat module),
${}_B\Sigma_A={}_BA^{(n)}_A$ for a positive integer $n$, and
${}_B\coring{D}$ and $\coring{D}_B$ are flat modules, then the
functors $-\tensor{B}A^{(n)}: \rcomod{\coring{D}} \rightarrow
\rcomod{A^{(n)}[\coring{D}]}$ and $A^{(n)}\tensor{B}-:
\lcomod{\coring{D}} \rightarrow \lcomod{A^{(n)}[\coring{D}]}$ are
equivalences  (i.e.\ $\coring{D}$ is Morita-Takeuchi equivalent to
$A^{(n)}[\coring{D}]$).
\end{corollary}
\begin{proof}
In the notation of Theorem~\ref{thm.equiv}, take $\cC \Sigma[\cD]$ and $\sigma$ the identity map.  Observe that since,
for all right $\cD$-comodules $M$,
$M\Box_\cD(\cD\otimes_B\Sigma)\simeq M\otimes_B \Sigma$, the fact
that $\Sigma$ is a faithfully flat left $B$-module implies that
$\cD\otimes_B\Sigma$ is a faithfully coflat left $\cD$-comodule.
Finally, in this case $-\ot_B\Sigma$ is the push-out functor and
$-\cotensor{\Sigma[\coring{D}]}(\Sigma^*\ot_B\cD)$ is the
pull-back functor, hence the assertions follow from
Theorem~\ref{thm.equiv}.
\end{proof}

\subsection{Generalised descent}

To any ring extension $\alpha: B \rightarrow A$ one can
associate a category of {\em right descent data} $\Desc_{\alpha}$
defined in \cite[3.3]{Nuss:1997} (cf.\ \cite{Cipolla:1976},
\cite{Knus/Ojanguren:1974}). As observed in \cite[Example 1.2]{Brz:str},
 the category $\Desc_{\alpha}$ is
isomorphic to the category of right comodules over the canonical
Sweedler $A$-coring $A\tensor{B}A$. This isomorphism allows one to formulate
a generalised
descent theorem \cite[Theorem 5.6]{Brz:str}. In this subsection we
introduce the category of descent data relative to a coring
morphism, and give a descent theorem in this general case.

Let $(\gamma,\alpha): (\cD\!:\!B) \rightarrow (\cC\!:\!A)$ be a coring
morphism, and consider the base ring extension coring
$A_{\alpha}[\cD]=A\ot_B\cD\ot_BA$.
There are two coring morphisms associated to $(\gamma,\alpha)$,
\[
\xymatrix@R=0pt@C=20pt{\widetilde{\gamma}:A_{\alpha}[\cD]
\ar@{->}[r] & \cC, \\ a\tensor{B}d\tensor{B}a' \ar@{|->}[r] &
a\gamma(d)a', } \xymatrix@R=0pt@C=20pt{
(\widetilde{\alpha},\alpha): (\cD\!:\!B) \ar@{->}[r] &
(A_{\alpha}[\cD]\!:\!A), \\ (d,b) \ar@{|->}[r] &
(1\tensor{B}d\tensor{B}1,\alpha(b)), }
\]
such that
\begin{equation}\label{cuatro}
\xymatrix@R=30pt@C=70pt{ \cD \ar@{->}^-{\widetilde{\alpha}}[r]
\ar@{->}_-{\gamma}[dr]
 & A_{\alpha}[\cD] \ar@{->}^-{\eps_{\cD,A_\alpha}}[r]
 \ar@{->}^-{\widetilde{\gamma}}[d] & A\ot_{B}A \\
& \cC & }
\end{equation}
is a commutative diagram of coring morphisms (cf.\ \cite[Section~24]{BrzWis:cor}).

A {\em descent datum} associated to a coring morphism $(\gamma,\alpha): (\cD\! :\! B) \rightarrow (\cC\! :\! A)$ is a pair $(X,\rho_X')$ consisting
of a right $A$-module $X$ and a right $A$-linear map $\rho_X': X
\rightarrow X\tensor{B}\coring{D}\tensor{B}A$ (here $X$ is considered
as a right $B$-module by restriction of scalars) rendering commutative the following
diagrams
\[
\xymatrix@R=30pt@C=70pt{ X \ar@{->}^-{X}[rr]
\ar@{->}^-{\rho_X'}[d] & &
X  \\
X\tensor{B}\coring{D}\tensor{B}A
\ar@{->}^-{X\tensor{B}\eps_{\coring{D}}\tensor{B}A}[r] &
X\tensor{B}B\tensor{B}A \ar@{->}^-{\simeq}[r] & X\tensor{B}A
\ar@{->}^-{\sigma^l_X}[u]}
\]
and
\[
\xymatrix@R=30pt@C=70pt{ X \ar@{->}^-{\rho_X'}[rr]
\ar@{->}^-{\rho_X'}[d]& & X\tensor{B}\coring{D}\tensor{B}A
\ar@{->}|-{\rho_X'\tensor{B}\coring{D}\tensor{B}A}[d] \\
X\tensor{B}\coring{D}\tensor{B}A
\ar@{->}^-{X\tensor{B}\Delta_{\coring{D}}\tensor{B}A}[r] &
X\tensor{B}\coring{D}\tensor{B}\coring{D}\tensor{B}A
\ar@{->}^-{X\tensor{B}\nu_{\coring{D},\coring{D}}\tensor{B}A}[r] &
X\tensor{B}\coring{D}\tensor{B}A_{\alpha}[\cD]. }
\]
Here $\sigma_{-}^l: -\tensor{B}A \rightarrow -$ is the natural
transformation defined by right multiplication, and $\nu_{-,-}:
-\tensor{B}- \rightarrow -\tensor{B}A\tensor{B}-$ is the canonical
natural transformation. A {\em morphism} of descent data is a
right $A$-linear map $f: (X,\rho_X) \rightarrow (Z,\rho_Z')$ such
that
\[
\xymatrix@C=60pt@R=30pt{ X \ar@{->}^-{f}[r] \ar@{->}_-{\rho_X'}[d]
&
Z \ar@{->}^-{\rho'_Z}[d] \\
X\tensor{B}\coring{D}\tensor{B}A
\ar@{->}^-{f\tensor{B}\coring{D}\tensor{B}A}[r] &
Z\tensor{B}\coring{D}\tensor{B}A }
\]
is a commutative diagram. Descent data associated to a coring morphism $(\gamma,\alpha)$ and their morphisms
form a category called the category of {\em right descent data
associated to $(\gamma,\alpha)$} and denoted by
$\Desc_{(\gamma,\alpha)}$.

\begin{lemma}\label{desc-comod}
Let $(\gamma,\alpha): (\cD\!:\!B) \rightarrow (\cC\!:\!A)$ be a
morphism of corings. Then $\Desc_{(\gamma,\alpha)}$ is isomorphic to  $\rcomod{A_{\alpha}[\cD]}$.
\end{lemma}
\begin{proof}
The  isomorphism of categories is provided by the following two
functors. Given any right $A_{\alpha}[\coring{D}]$-comodule $X$,
define a right $A$-linear map
\[
\xymatrix{\rho_X': X \ar@{->}^-{\varrho^X}[r] &
X\tensor{A}A\tensor{B}\coring{D}\tensor{B}A \ar@{->}^-{\simeq}[r] &
X\tensor{B}\coring{D}\tensor{B}A .}
\]
It is clear that $(X,\rho'_X)$ is an object of
$\Desc_{(\gamma,\alpha)}$. Obviously, any right
$A_{\alpha}[\cD]$-colinear map induces a morphism in
$\Desc_{(\gamma,\alpha)}$. Conversely, let $(X,\rho_X')$ be an
object of $\Desc_{(\gamma,\alpha)}$. Then $X$ is a right
$A_{\alpha}[\cD]$-comodule with the coaction
\[
\xymatrix{\varrho^X: X \ar@{->}^-{\rho_X'}[r] &
X\tensor{B}\coring{D}\tensor{B}A \ar@{->}^-{\simeq}[r] &
X\tensor{A}A\tensor{B}\coring{D}\tensor{B}A. }
\]
Clearly any arrow of $\Desc_{(\gamma,\alpha)}$ induces a right
$A_{\alpha}[\cD]$-colinear map. Finally, a straightforward computation shows that
the constructed  functors are mutually inverse.
\end{proof}

\begin{corollary}
Let $(\gamma,\alpha) :(\cD\!:\!B) \rightarrow (\cC\!:\!A)$ be a coring
morphism. If ${}_B\coring{D}$ is a flat module and ${}_BA$ is a
faithfully flat module, then
\[
\xymatrix@C=70pt{ \rcomod{\coring{D}} \ar@{->}^-{-\tensor{B}A}[r]
& \rcomod{A_{\alpha}[\cD]} \equiv \Desc_{(\gamma,\alpha)}}
\]
is an equivalence of categories.
\end{corollary}
\begin{proof}
This corollary is a straightforward consequence of Lemma~\ref{desc-comod} and Corollary~\ref{cor.equiv.com}.
\end{proof}

The following commutative diagram of functors summarises and
combines the old and the new situations
\[
\xymatrix@R=40pt{ \rcomod{\coring{D}} \ar@{->}^-{-\tensor{B}A}[rr]
\ar@{->}_-{U_B}[dd] \ar@{->}^-{-\tensor{B}A}[rd] & &
\rcomod{\coring{C}} \ar@{->}^-{U_A}[dd]
\\ &
\rcomod{A_{\alpha}[\cD]} \simeq \Desc_{(\gamma,\alpha)}
\ar@{->}^-{(-)_{\widetilde{\gamma}}}[ru]
\ar@{->}^<<<<<{(-)_{\eps_{\cD,A_\alpha}}}[dd] &
\\ \rmod{B}
\ar@{->}^<<<<<<<<<<<<<{-\tensor{B}A}[rr] \ar@{->}_-{-\tensor{B}A}[rd] & & \rmod{A} \\
&\rcomod{A\tensor{B}A} \simeq \Desc_{\alpha} \ar@{->}_-{u_A}[ru] &
}
\]
Here $-\tensor{B}A: \rcomod{\coring{D}} \rightarrow
\rcomod{\coring{C}}$ is the right push-out functor induced by the
module-morphism $(A_\alpha,\widetilde{\gamma})$, and
$(-)_{\widetilde{\gamma}}$ and  $(-)_{\eps_{\cD,A_\alpha}}$ are the
induction functors associated to the coring morphisms
$\widetilde{\gamma}$ and $\eps_{\cD,A_\alpha}$, respectively (the latter
defined in diagram~\eqref{cuatro}).

\section*{Acknowledgements}
T.\ Brzezi\'nski started working on this paper during his visit to
Santiago de Compostela. He would like to thank Emilio Villanueva
Novoa and other members of the Department of Algebra of the
University of Santiago de Compostela for very warm hospitality.

This research has been partly developed during the postdoctoral
stay of L. El Kaoutit at the Universidad de Granada supported by
the grant SB2003-0064 from the Mi\-nis\-terio de Educaci{\'o}n,
Cultura y Deporte of Spain.

 This work is also partially supported by the
grant MTM2004-01406 from the Ministerio de Educaci\'{o}n y Ciencia of
Spain.

All the authors would like to thank Gabriella B\"ohm for discussions.

  \end{document}